\documentclass[12pt]{article} 
\usepackage{amssymb, amsmath, amsthm}
\usepackage{geometry} 
\usepackage{amsfonts} 
\usepackage{amsmath} 
\usepackage[usenames, dvipsnames]{color}
\usepackage{relsize}
\usepackage{graphicx}
\usepackage[utf8]{inputenc}

    \newcommand{\rn}{\R^n}
    \newcommand{\K}{{\mathcal K}}

\newtheorem{theorem}{Theorem}[section]
\newtheorem{lemma}[theorem]{Lemma}
\newtheorem{remark}[theorem]{Remark}
\newtheorem{proposition}[theorem]{Proposition}
\newtheorem{definition}[theorem]{Definition}

    \newcommand{\sub}{\subseteq}

    \newcommand{\ct}[1]{\langle {#1}\rangle \lower.3ex\hbox{$_{t}$}}
    \newcommand{\lt}[1]{[ {#1}] \lower.3ex\hbox{$_{t}$}}

\newcommand{\cn}{\mathcal{C}^n}

\newcommand{\R}{{\bf R}}
\newcommand{\N}{{\bf N}}

\newcommand{\C}{{\mathcal C}}

\newcommand{\ml}{{\mathcal H}^{n}}

\newcommand{\cl}{{\rm cl}}
\newcommand{\interno}{{\rm int}}
\newcommand{\clos}{{\rm cl}}
\newcommand{\dom}{{\rm dom}}
\newcommand{\dimension}{{\rm dim}}
\newcommand{\relint}{{\rm relint}}

\newcommand{\epi}{{\rm epi}}
\newcommand\restr[2]{{
  \left.\kern-\nulldelimiterspace 
  #1 
  \right|_{#2} 
  }}
\newcommand\ext[2]{{
  \left.\kern-\nulldelimiterspace 
  #1 
  \right|^{#2} 
  }}

\newtheorem{corollary}[theorem]{Corollary}

\begin{document}


\title{Monotone valuations\\ on the space of convex functions}



\author{L. Cavallina and A. Colesanti\footnote{Supported by G.N.A.M.P.A. and by the FIR project 2013 
``Geometric and Qualitative aspects of PDE's'' }}


\date{}
\maketitle

\begin{abstract} We consider the space $\C^n$ of convex functions $u$ defined in $\R^n$ with values in $\R\cup\{\infty\}$, which are lower
semi-continuous and such that $\lim_{|x|\to\infty} u(x)=\infty$. We study the valuations defined on $\C^n$ which are invariant under the composition 
with rigid motions, monotone and verify a certain type of continuity. Among these valuations we prove integral representation formulas for those which are,
additionally, simple or homogeneous.   
\end{abstract}
\bigskip


\bigskip


\tableofcontents


\section{Introduction} 

The aim of this paper is to begin an exploration of the valuations defined on the space of convex functions, having as a model the valuations 
of convex bodies.

\medskip

The theory of valuations is currently a significant part of convex geometry. 
We recall that, if $\K^n$ denotes the set of convex bodies (compact and convex sets) in $\R^n$, a (real-valued) valuation is an application $\sigma\,:\,\K^n\to\R$ 
that verifies the following (restricted) additivity condition
\begin{equation}\label{intro1}
\sigma(K\cup L)+\sigma(K\cap L)=\sigma(K)+\sigma(L)\quad\forall\,K,L\in\K^n\;\mbox{such that $K\cup L\in\K^n$,}
\end{equation}
together with
\begin{equation}\label{intro2}
\sigma(\emptyset)=0.
\end{equation}
 
In the realm of convex geometry the most familiar examples of valuations are the so-called {\em intrinsic volumes} $V_k$, $k\in\{0,1,\dots,n\}$,
which have many 
additional properties such as: invariance under rigid motions, continuity with respect to the Hausdorff metric, homogeneity and monotonicity. 
Note that intrinsic volumes include the volume (here denoted by $V_n$) itself, i.e. the Lebesgue measure, which is clearly a valuation. 

A celebrated result by Hadwiger (see \cite{Hadwiger}, \cite{Klain}, \cite{Klain-Rota}, \cite{Schneider}) provides a characterization of an important class
of valuations on $\K^n$.

\begin{theorem}[{\bf Hadwiger}] Every rigid motion invariant valuation on $\K^n$, which is continuous
(with respect to the Hausdorff metric) or monotone, can be written as the linear combination of intrinsic volumes. 
\end{theorem}

A special case of this theorem is known as the {\em volume theorem}.

\begin{corollary}\label{volt} Let $\sigma$ be a rigid motion invariant and continuous (or monotone) valuation on $\K^n$, which is simple, i.e. 
$\sigma(K)=0$ for every $K$ such that $\dim(K)< n$. Then $\sigma$ is a multiple of the volume: there exists a constant $c\in\R$ such that
$$
\sigma(K)=c V_n(K)\quad\forall\, K\in\K^n.
$$
\end{corollary}

These deep results gave a strong impulse to the development of the theory which, in the last decades, was enriched by a wide variety of new
results and counts now a considerable number of prolific ramifications. A survey on the state of the art of this subject is presented in the monograph \cite{Schneider} 
by Schneider (see chapter 6), along with a detailed list of references. 

\medskip

Recently the study of valuations was extended from spaces of sets, like $\K^n$, to spaces of functions. The condition \eqref{intro1} is adapted
to this situation replacing union and intersection by ``max'' and ``min''. In other words, if $\mathcal X$ is a space of functions, an application 
$\mu\,:\,{\mathcal X}\to\R$ is called a valuation if 
\begin{equation}\label{intro2}
\mu(u\vee v)+\mu(u\wedge v)=\mu(u)+\mu(v)\quad
\end{equation}
for every $u,v\in{\mathcal X}$ such that $u\wedge v\in{\mathcal X}$ and $u\vee v\in{\mathcal X}$.
Here $u\vee v$ and $u\wedge v$ denote the point-wise maximum and minimum of $u$ and $v$, respectively. To motivate this definition one can 
observe that the epigraphs of $u\vee v$ and $u\wedge v$ are the intersection and the union of the epigraphs of $u$ and $v$,
respectively. The same property is shared by sub-level sets: for every $t\in\R$ we have
$$
\{u\vee v<t\}=\{u<t\}\cap\{v<t\}\quad\mbox{and}\quad
\{u\wedge v<t\}=\{u<t\}\cup\{v<t\}.
$$
This fact will be particularly important in the present paper.  

Valuations defined on the Lebesgue spaces $L^p(\R^n)$, $p\ge1$, were studied by Tsang in \cite{Tsang-2010}. In relation
to some of the results presented in our paper it is interesting to mention that one of the results of Tsang asserts that any translation invariant 
and continuous valuation $\mu$ on $L^p(\R^n)$ can be written in the form
\begin{equation}\label{intro3}
\mu(u)=\int_{\R^n} f(u(x))\,dx\quad\forall\, u\in L^p(\R^n)
\end{equation}
where $f$ is a continuous function subject to a suitable growth condition at infinity. The results of Tsang have been extended to
Orlicz spaces by Kone, in \cite{Kone}. 
Valuations of different types (taking values in $\K^n$ or in spaces of matrices, instead of $\R$),
defined on Lebesgue, Sobolev and BV spaces, have been considered in \cite{Tsang-2011}, 
\cite{Ludwig-2011a}, \cite{Ludwig-2012}, \cite{Ludwig-2013}, \cite{Wang-thesis}, \cite{Wang} and \cite{Ober}
(see also \cite{Ludwig-2011} for a survey).

\medskip

Wright, in his PhD Thesis \cite{Wright} and subsequently in collaboration with Baryshnikov and Ghrist in \cite{BGW}, studied a rather different class,
formed by the so-called {\em definable functions}. We cannot give here the details of the construction of these functions, but we mention that the main result
of these works is a characterization of valuation as suitable integrals of intrinsic volumes of level sets. These type of integrals 
will have a crucial role in our paper, too.  

\medskip

The class of functions that we are considering is 
$$
\C^n=\{u\,:\,{\bf R}^n\to{\bf R}\cup\{\infty\},\,\mbox{$u$ convex, l.s.c.,
$\lim_{|x|\to\infty}u(x)=\infty$}\}.
$$
It includes the so-called indicatrix functions of convex bodies, i.e. functions of the form
$$
I_K\,:\,\R^n\to\R\cup\{\infty\},\quad
I_K(x)=\left\{
\begin{array}{ll}
\mbox{$0$ if $x\in K$,}\\
\mbox{$\infty$ if $x\notin K$,}
\end{array}
\right.
$$ 
where $K$ is a convex body. Note that the function $\boldsymbol{\infty}$, identically equal to $\infty$, belongs to $\C^n$. This element will play 
in some sense the role of the empty set. If $u\in\C^n$ we will denote by $\dom(u)$ the set where $u$ is finite; 
if $u\ne\boldsymbol{\infty}$, this is a non-empty convex set, and then its dimension $\dim(\dom(u))$ is well defined and it is an integer between $1$ and 
$n$.

We say that a functional $\mu\,:\,\C^n\to\R$ is a valuation if it verifies \eqref{intro2} for every $u,v\in\C^n$ such that
$u\wedge v\in\C^n$ (note that $\C^n$ is closed under ``$\vee$''), and $\mu(\boldsymbol{\infty})=0$. We are interested in 
valuations which are rigid motion invariant (i.e. $\mu(u)=\mu(u\circ T)$ for every $u\in\C^n$ and every rigid motion $T$ of 
$\R^n$), and monotone decreasing (i.e. $\mu(u)\le\mu(v)$ for every $u,v\in\R^n$ such that $u\ge v$ point-wise in $\R^n$). 
As $\mu(\boldsymbol{\infty})=0$, we immediately have that they are non-negative in $\C^n$. We will also need to introduce a notion 
of continuity of valuations. In this regard, note that for (rigid motion invariant) valuations on the space of convex bodies, continuity and monotonicity are 
conditions very close to each other. In particular monotonicity implies continuity, as Hadwiger's theorem shows. The situation on $\C^n$ is rather 
different and it is easy to provide examples of monotone valuations which are not continuous with respect to any reasonable notion of convergence on
$\C^n$. We say that a valuation on $\C^n$ is {\em monotone-continuous} ({\em m-continuous} for short) if 
$$
\lim_{i\to\infty}\mu(u_i)=\mu(u)
$$
whenever $u_i$, $i\in\N$, is a decreasing sequence in $\C^n$ converging to $u\in\C^n$ point-wise in the relative interior of
$\dom(u)$, and such that $u_i\ge u$ in $\R^n$ for every $i$.  

\bigskip

How does a ``typical'' valuation of this kind look like? A first answer is provided by functionals of type \eqref{intro2}; indeed we will see in section
\ref{sec integral valuations} that if 
$f\,:\,\R\to\R$ is a non-negative decreasing function which verifies the integrability condition
\begin{equation}\label{intro4}
\int_0^\infty f(t)t^{n-1}dt<\infty
\end{equation}
then the functional
\begin{equation}\label{intro5}
\mu(u)=\int_{\dom(u)}f(u(x))dx
\end{equation}
is a rigid motion invariant, monotone decreasing valuation on $\C^n$, which is moreover m-continuous if $f$ is right-continuous. A valuation 
of this form vanishes obviously on every function $u\in\C^n$ such that $\dim(\dom(u))<n$: 
\begin{equation}\label{intro6}
\dim(\dom(u))<n\quad\Rightarrow\quad \mu(u)=0.
\end{equation}
When $\mu$ has this property we will say that it is {\em simple}. Our first characterization result is the following theorem, proven in section 
\ref{Characterization results I: simple valuations}.

\begin{theorem}\label{thm intro 1}
Let $\mu$ be a valuation on $\C^n$ which is: rigid motion invariant, decreasing, m-continuous and simple. Then there exists a decreasing function 
$f$ defined on $\R$, left-continuous and verifying \eqref{intro4}, such that $\mu$ can be written in the form \eqref{intro5}.
\end{theorem}

The proof of this fact is based on a rather simple idea, even if there are several technical points to transform it into a rigorous argument.
First, $\mu$ determines the function $f$ as follows: for $t\in\R$ let $\sigma_t\,:\,\K^n\to\R$ be defined as
$$
\sigma_t(K)=\mu(t+I_K).
$$
It is straightforward to check that this is a rigid motion and monotone increasing valuation, so that by the volume theorem there exists a constant, which will
depend on $t$ and which we call $f(t)$, such that
$$
\mu(t+I_k)=\sigma_t(K)=f(t) V_n(K)\quad\forall\,K\in\K^n,\,
\forall\,t\in\R.
$$
As $\mu$ is decreasing, $f$ is decreasing too. Now for every $u\in\C^n$ and $K\sub\dom(u)$, by monotonicity we obtain
$$
f(\max_K u) V_n(K)=\mu(\max_K u+I_K)\le
\mu(u+I_k)\le f(\min_K u) V_n(K)=\mu(\min_K u+I_K).
$$
This chain of inequalities, and the fact that $\mu$ is simple, permit to compare easily the value of $\mu(u)$ 
with upper and lower Riemann sums of $f\circ u$, over suitable partitions of subsets of $\dom(u)$. This leads to the proof of 
\eqref{intro5}. Monotonicity is an essential ingredient of this argument. It would be very interesting to obtain a 
similar characterization of simple valuations without this assumption.

\medskip

If we apply the layer cake (or Cavalieri) principle, we obtain a second way 
of writing the valuation $\mu$ defined by \eqref{intro5}:
\begin{equation}\label{intro7}
\mu(u)=\int_{\R}V_n(\cl(\{u<t\}))d\nu(t)
\end{equation}
where $V_n$ is the $n$-dimensional volume, ``$\cl$'' is the closure, and $\nu$ is a Radon measure on $\R$ identified by the equality
$$
f(t)=\nu((t,\infty))\quad\forall\, t\in\R.
$$ 
Note that the set $\cl(\{u<t\})$ is a compact convex set, i.e. a convex body. Formula \eqref{intro7} suggests to consider the more
general expression
\begin{equation}\label{intro8}
\mu(u)=\int_{\R}V_k(\cl(\{u<t\}))d\nu(t)
\end{equation}
where, for $k\in\{0,\dots,n\}$, $V_k$ is the $k$-th intrinsic volume and $\nu$ is a Radon measure on $\R$. We will see 
(still in section \ref{sec integral valuations}) that the integral in \eqref{intro8} 
is finite for every $u\in\C^n$ if and only if
\begin{equation}\label{intro9}
\int_0^\infty t^kd\nu(t)<\infty
\end{equation}
(which is equivalent to \eqref{intro4} when $k=n$) and in this case it defines a rigid motion invariant, decreasing and m-continuous valuation
on $\C^n$.  Valuations of type \eqref{intro8} are homogeneous of order $k$ in the following sense. For $u\in\C^n$ and $\lambda>0$, let
$u_\lambda\,:\,\R^n\to\R\cup\{\infty\}$ be defined by
$$
u_\lambda(x)=u\left(\frac x\lambda\right)\quad\forall\, x\in\R^n
$$
(note that $u_\lambda\in\C^n$). Then
$$
\mu(u_\lambda)=\lambda^k\mu(u).
$$
In section \ref{Characterization results II: homogeneous valuations} we will prove the following fact.

\begin{theorem}\label{thm intro 2} 
Let $\mu$ be a valuation on $\C^n$ wich is rigid motion invariant, decreasing, m-continuous and $k$-homogeneous for some 
$k$. Then necessarily $k\in\{0,1,\dots,n\}$ and there exists a Radon measure $\nu$ on $\R$, verifying \eqref{intro9}, such that
$\mu$ can be written in the form \eqref{intro8}.
\end{theorem}

\medskip

Theorems \ref{thm intro 1} and \ref{thm intro 2} may suggest that valuations of type \eqref{intro8} could form a sort of generators for 
invariant, monotone and m-continuous valuations on $\C^n$, playing a similar role to intrinsic volumes for convex bodies (with the difference that 
the dimension of the space of valuations on $\C^n$ is infinite). On the other hand in the conclusive section of the paper we show the existence of
valuations on $\C^n$ with the above properties, which cannot be decomposed as the sum of  homogeneous valuations, and hence are not the sum of valuations of the form \eqref{intro8}.

\medskip

The second author would like to thank Monika Ludwig for encouraging his interest towards valuations on convex functions 
and for the numerous and precious conversations that he had with her on this subject.

\section{Preliminaries}\label{sec-preliminaries}

We work in the $n$-dimensional Euclidean space ${\bf R}^n$, $n\ge1$, endowed with the usual Euclidean norm $|\cdot|$ and 
scalar product $(\cdot,\cdot)$. For $x_0\in\R^n$ and $r>0$, $B_r(x_0)$ denotes the closed ball centred at $x_0$ with radius $r$;
when $x_0=0$ we simply write $B_r$.
For $k\in[0,n]$, the $k$-dimensional Hausdorff measure is denoted by $\mathcal{H}^k$. In particular $\ml$ denotes the Lebesgue measure 
in ${\bf R}^n$ (which, as we said, will be often indicated by $V_n$, especially when referred to convex bodies).
Integration with respect to such measure will be always denoted simply by $dx$, where $x$ is the integration variable.
Given a subset $A$ of $\R^n$ we denote by $\interno(A)$ and $\clos(A)$ its interior and its closure, respectively.

As usual, we will denote by ${\bf O}(n)$ and ${\bf SO}(n)$ respectively, the group of rotations  and of proper rotations of $\R^n$. By a {\em rigid motion} we mean 
the composition of a rotation and a translation, i.e. a mapping $T\,:\,\R^n\to\R^n$ such that there exist $R\in{\bf O}(n)$ and $x_0\in\R^n$ for which
$$
T(x)=R(x)+x_0,\quad\forall\, x\in\R^n.
$$

\subsection{Convex bodies}\label{subsection convex bodies}

A {\em convex body} is a compact convex subset of ${\bf R}^n$. We will denote by $\K^n$ the family of convex bodies in ${\bf R}^n$. For all 
the notions and results concerning convex bodies we refer to the monograph \cite{Schneider}. The set $\K^n$ can be endowed with a metric,
induced by the Hausdorff distance (see \cite{Schneider} for the definition).


Let $K\in\K^n$; if $\interno(K)=\emptyset$, then $K$ is contained in some $k$-dimensional affine sub-space
of $\R^n$, with $k<n$; the smallest $k$ for which this is possible is called the dimension of $K$, and is denoted by $\dimension(K)$. Clearly,
if $K$ has non-empty interior we set $\dimension(K)=n$. Using this notion we can define the relative interior of $K$ as the subset of those points
$x$ of $K$ for which there exists a $k$-dimensional ball centred at $x$ and contained in $K$, where $k=\dimension(K)$. The relative 
interior will be denoted by $\relint(K)$. The notion of relative interior can be given in the same way for every convex subset of $\R^n$.

\medskip

$\K^n$ can be naturally equipped with an addition (Minkowski, or vector, addition) and a multiplication by non-negative reals. Given $K,L\in\K^n$
and $s\ge0$ we set
$$
K+L=\{x+y\,:\,x\in K,\quad y\in L\}
$$
and
$$
sK=\{sx\,:\,x\in K\}.
$$
$\K^n$ is closed with respect to these operations. 

\medskip

To every convex body $K\in\K^n$ we may assign a sequence of $(n+1)$ numbers, $V_k(K)$, $k=0,\dots,n$, 
called the {\it intrinsic volumes} of $K$; for their definition see \cite[Chapter 4]{Schneider}. We recall in 
particular that $V_n(K)$ is the volume, i.e. the Lebesgue measure, of $K$, while $V_0(K)=1$ for every $K\in\K^n\setminus\{\emptyset\}$.
More generally, if $K$ is a convex body in
$\R^n$ having dimension $k\in\{0,1,\dots,n\}$, then $V_k(K)$ is the $k$-dimensional Lebesgue measure of $K$
as a subset of $\R^k$. 
As real-valued functionals defined on $\K^n$, intrinsic volumes are continuous, monotone increasing with respect to 
set inclusion and invariant under the action of rigid motions: $V_i(T(K))=V_i(K)$ for every $i\in\{0,\dots,n\}$, $K\in\K^n$ and
for every rigid motion $T$. Moreover, the intrinsic volumes are special and important examples of 
{\em valuations} on the space of convex bodies. We recall that a (real-valued) valuation on $\K^n$ is a mapping $\sigma\,:\,\K^n\to\R$
such that $\sigma(\emptyset)=0$ and 
$$
\sigma(K\cup L)+\sigma(K\cap L)=\sigma(K)+\sigma(L)\quad
\mbox{$\forall K,L\in\K^n$ such that $K\cup L\in\K^n$.}
$$

The following characterization theorem of Hadwiger (see, for instance, \cite[Chapter 6]{Schneider}) will be a crucial tool in this paper.

\begin{theorem}\label{Hadwiger-bodies} Let $\sigma$ be a valuation on $\K^n$ which is invariant with respect to rigid motions, and
either continuous with respect to the Hausdorff metric or monotone with respect to set inclusion. Then $\sigma$ is the linear combination 
of intrinsic volumes, i.e. there exist $c_0,\dots,c_n\in\R$, such that
$$
\sigma(K)=\sum_{i=0}^n c_i V_i(K)\quad\mbox{for every $K\in\K^n$.}
$$
Moreover, if $\sigma$ is increasing (resp. decreasing) then $c_i\ge0$ (resp. $c_i\le0$) for every $i\in\{0,\dots,n\}$.
\end{theorem}

\section{The space $\C^n$}

Let us a consider a function $u\,:\,\R^n\to\R\cup\{\infty\}$, which is {\em convex}. We denote the so-called {\it domain} of $u$ as
$$
\dom(u)=\{x\in\R^n\,:\, u(x)<\infty\}.
$$
By the convexity of $u$, $\dom(u)$ is a convex set. By standard properties
of convex functions, $u$ is continuous in the interior of $\dom(u)$
and it is Lipschitz continuous in any compact subset of $\interno(\dom(u))$. We will sometimes use 
the following notation, for $u\in\C^n$,
$$
\Omega_u=\interno(\dom(u)).
$$

\medskip

In this work we focus in particular on the following space of convex functions:
\begin{equation}\label{I.1}
\C^n=\{u\,:\,{\bf R}^n\to{\bf R}\cup\{\infty\},\,\mbox{$u$ convex, l.s.c.,
$\lim_{|x|\to\infty}u(x)=\infty$}\}.
\end{equation}
Here by l.s.c. we mean {\em lower semi-continuous}, i.e.
$$
\liminf_{x\to x_0} u(x)\ge u(x_0)\quad\forall
x_0\in\R^n.
$$
Note that the function $\boldsymbol{\infty}$ (which, we recall, is identically equal to $\infty$ on $\rn$) belongs to our functions space. As it will be clear in the sequel, this special function plays the 
role that the empty set has for valuations defined on families of sets (instead of functions). 

\begin{remark} {\em Let $u\in\C^n$. As a consequence of convexity and the behavior at infinity we have that
$$
\inf_{\R^n}u>-\infty.
$$
Moreover, by the lower semi-continuity, $u$ admits a minimum in $\R^n$. We will often use the notation
$$
m(u)=\min_{\R^n}u\,.
$$
}\end{remark}

We will also need to consider the following subset of $\C^n$:
$$
\C^n_b=\{u\in\C^n\,:\,\mbox{$\dom(u)$ is bounded}\}.
$$

\bigskip

Let $A\sub\R^n$; we denote by $I_A\,:\,\R^n\to\R\cup\{\infty\}$ the so-called {\em indicatrix function} of $A$, which is defined by
$$
I_A=\left\{
\begin{array}{ll}
0\quad&\mbox{if $x\in A$,}\\
\infty\quad&\mbox{if $x\notin A$.} 
\end{array}
\right.
$$
If $K\subset\R^n$ is a convex body, then $I_K\in\C^n$.

\bigskip




Sub-level sets of functions belonging to $\C^n$ will be of fundamental importance in this paper. Given $u\in\C^n$ and $t\in\R$ we set
$$
K_t:=\{u\le t\}=\{x\in\R^n\,:\,u(x)\le t\},
$$
and
$$
\Omega_t:=\{u<t\}=\{x\in\R^n\,:\,u(x)<t\}.
$$
Both sets are empty for $t<m(u)$. $K_t$ is a convex body for all $t\in\R$, by the properties of $u$. 
For all real $t$, $\Omega_t$ is a bounded (possibly empty)
convex set, so that its closure $\cl(\Omega_t)$ is a convex body, obviously contained in $K_t$.

\begin{lemma}\label{lemma sottolivelli}
Let $u\in\C^n$; for every $t>m(u)$ 
$$
\relint(K_t)\sub\Omega_t.
$$
\end{lemma}

\begin{proof} We start by considering the case in which $\dim(K_t)=n$. Assume by contradiction that there exists a point 
$x\in\interno(K_t)$ such that $u(x)=t$. Then $x$ is a local maximum for $u$ but, by convexity, this is possible only if $u\equiv t$ in $K_t$,
which, in turn implies that $t=m(u)$, a contradiction. 

If $\dim(K_t)=k<n$ then, by convexity, $\dom(u)$ is contained in a $k$-dimensional affine subspace $H$ of $\R^n$, and we can apply the previous argument
to $u$ restricted to $H$ to deduce the assert of the lemma.
\end{proof}

\begin{corollary}\label{corollario sottolivelli} Let $u\in\C^n$; for every $t>m(u)$ 
$$
\cl(\Omega_t)=K_t. 
$$
\end{corollary}

\subsection{On the intrinsic volumes of sub-level sets}\label{On the intrinsic volumes of sub-level sets}
As we have just seen, if $u\in\C^n$ and $t\in\R$, the set
$$
\Omega_t=\{u<t\}
$$
is empty for $t\le m(u)$ and it is a bounded convex set for $t> m(u)$. For $k\in\{0,\dots,n\}$, we define the function 
$v_k(u;\cdot)\,:\,\R\to\R$ as follows
$$
v_k(u;t)=V_k(\cl(\Omega_t)).
$$
As intrinsic volumes are non-negative and monotone with respect to set inclusion and 
the set $\Omega_t$ is increasing with respect to inclusion as $t$ increases, $v_k(u;\cdot)$ is a non-negative increasing function. In particular it is a function 
of bounded variation, so that there exists a (non-negative) Radon measure on $\R$, that we will denote by $\beta_k(u;\cdot)$, which represents
the weak, or distributional, derivative of $v_k$ (see for instance \cite{AFP}). 

We want to describe in a more detailed way the structure of the measure $\beta_k$. In general, the measure representing the weak derivative 
of a non-decreasing function consists of three parts: a jump part, a Cantor like part and an absolutely continuous part (with respect to Lebesgue measure).
We will see that $\beta_k$ does not have a Cantor part and its jump part, if any, is a single Dirac delta at $m(u)$.

As a starting point, note that as $v_i$ is identically zero in $(-\infty,m(u)]$ then $\beta_k(u;\eta)=0$ for every measurable set 
$\eta\sub(-\infty,m(u))$. 
On the other hand, in $(m(u),\infty)$, due to the Brunn-Minkowski inequality for intrinsic volumes, the function $v_k$ have a more regular behavior  
than that of a non-decreasing function. Indeed, for $k\ge1$, let $t_0,t_1\in(m(u),\infty)$ and consider, for $\lambda\in[0,1]$, $t_\lambda=(1-\lambda)t_0+\lambda t_1$.  
Then we have the set inclusion
$$
K_{t_\lambda}\supseteq(1-\lambda)K_{t_0}+\lambda K_{t_1},
$$
which follows from the convexity of $u$. 
By the monotonicity of intrinsic volumes and the Brunn-Minkowski inequality for such functionals (see 
\cite[Chapter 7]{Schneider}), and by Corollary \ref{corollario sottolivelli}, we have
\begin{eqnarray*}
v_k(u;t_\lambda)&=&V_k(K_{t_\lambda})\ge
V_k((1-\lambda)K_{t_0}+\lambda K_{t_1})\\
&\ge&[(1-\lambda)V_k(K_{t_0})^{1/k}+\lambda V_k(K_{t_1})^{1/k}]^k\\
&=&[(1-\lambda)v_k(u;t_0)^{1/k}+\lambda v_k(u;t_1)^{1/k}]^k.
\end{eqnarray*}

In other words, the function $v_k$ to the power $1/k$ is concave in $(m(u),\infty)$. This implies in particular that $v_k$ is absolutely continuous
in $(m(u),\infty)$ so that the measure $\beta_k$ is absolutely continuous with respect to the Lebesgue measure in this interval, and its density is 
given by the point-wise derivative of $v_k$, which exists a.e. (see \cite[Chapter 3]{AFP}). 
Next we examine the behavior at $m(u)$; as $v_k$ is constantly zero in $(-\infty,m(u)]$
$$
\lim_{t\to m(u)^-}v_k(u;t)=0
$$ 
(in particular $v_k$ is left-continuous at $m(u)$). 
On the other hand let $t_i$, $i\in\N$, be a decreasing sequence converging to $m(u)$, with $t_i>m(u)$ for every $i$, 
and consider the corresponding sequence of convex bodies 
$L_i=\cl(\Omega_{t_i})$, $i\in\N$. This is a decreasing sequence and $L_i\supseteq K_{m(u)}$ for every $i$. Moreover, trivially
$$
K_{m(u)}=\bigcap_{i\in\N} L_i.
$$
This implies in particular that $K_{m(u)}$ is the limit of the sequence $L_i$ with respect to the Hausdorff metric (see \cite[Section 1.8]{Schneider}). Then
$$
\lim_{i\to\infty} V_k(L_i)=V_k(K_{m(u)})
$$
so that
$$
\lim_{t\to m(u)^+} v_k(u;t)=V_k(K_{m(u)})=V_k(\{u=m(u)\}).
$$
If
$$
V_k(K_{m(u)})>0
$$ 
then $v_k(u;t)$ has a jump discontinuity at $m(u)$ of amplitude $V_k(K_{m(u})$. In other words
$$
\beta_k(u;\{m(u)\})=V_k(K_{m(u)}).
$$ 
The case $k=0$ can be treated as follows: as $V_0(K)$ is the Euler characteristic of
$K$ for every $K$, i.e. is constantly 1 on $\K^n\setminus \{\emptyset\}$, $v_0(u;t)$ equals 0 for $t\le m(u)$ and equals 
1 for $t > m(u)$; hence $\beta_0$ is just the Dirac point mass measure concentrated at $m(u)$.

The following statement collects the facts that we have proven so far in this part.

\begin{proposition}\label{structure beta_i} 
Let $u\in\C^n$ and $k\in\{0,\dots,n\}$; let $K_t$, $v_k$ be defined as before. Define the measure $\beta_k$ as
$$
\beta_k(u;\cdot)=V_k(K_{m(u)})\delta_{m(u)}(\cdot)+\frac{dv_k}{dt}{\mathcal H}^1(\cdot),
$$
where $\delta$ denote the Dirac point-mass measure (and ${\mathcal H}^1$ is the Lebesgue measure on $\R$). Then $\beta_k(u;\cdot)$ is the
distributional derivative of $v_k$, more precisely
$$
v_k(u; t)=\beta_k(u;(m(u),t])\quad\forall\,t\ge m(u)\quad
\mbox{and $v_k(u; t)=0$ $\forall t\le m(u)$.}
$$
In particular $v_k(u; \cdot)$ is left-continuous at $m(u)$. 

\end{proposition}

\subsection{Max and min operations in $\C^n$}

As we will see, the definition of valuations on $\C^n$ is based on the point-wise minimum and maximum of convex functions. This part
is devoted to some basic properties of these operations.

Given $u$ and $v$ in $\C^n$ we set, for $x\in\R^n$,
$$
(u\vee v)(x)=\max\{u(x), v(x)\}=u(x)\vee v(x),\quad
(u\wedge v)(x)=\max\{u(x), v(x)\}=u(x)\wedge v(x).
$$
Hence $u\vee v$ and $u\wedge v$ are functions defined in $\R^n$, with values in $\R\cup\{\infty\}$.

\begin{remark}
{\em If $u,v\in\C^n$ then $u\vee v$ belongs to $\C^n$ as well. Indeed convexity and 
behavior at infinity  are straightforward. Concerning lower semicontinuity of $u\vee v$, this is equivalent to say 
that $\{u\vee v\le t\}$ is closed for every $t\in\R$, which follows immediately from the equality
$$
\{u\vee v\le t\}=\{u\le t\}\cup\{v\le t\}.
$$
On the contrary, $u,v\in\C^n$ does not imply, in general, that $u\wedge v\in\C^n$ (a counterexample is
given by the indicatrix functions of two disjoint convex bodies).} 
\end{remark}

\bigskip

For $u,v\in\C^n$ and $t\in\R$ the following relations are straightforward:
\begin{equation}\label{min-max 1}
\{u\le t\}\cap\{v\le t\}=\{u\vee v\le t\},\quad
\{u\le t\}\cup\{v\le t\}=\{u\wedge v\le t\};
\end{equation}
\begin{equation}\label{min-max 2}
\{u< t\}\cap\{v< t\}=\{u\vee v< t\},\quad
\{u< t\}\cup\{v< t\}=\{u\wedge v< t\}.
\end{equation}

In the sequel we will also need the following result.

\begin{proposition}\label{proposition min-max} Let $u,v\in\C^n$ be such that $u\wedge v\in\C^n$. Then, for every $t\in\R$,
\begin{eqnarray}
\cl(\{u< t\})\cap\cl(\{v< t\})&=&\cl(\{u\vee v< t\}),\label{min-max 3}\\
\cl(\{u< t\})\cup\cl(\{v< t\})&=&\cl(\{u\wedge v< t\}).\label{min-max 4}
\end{eqnarray}
\end{proposition}

\begin{proof} Equality \eqref{min-max 4} comes directly from the second equality in \eqref{min-max 2}, passing to the closures of the involved sets.
As for the proof of \eqref{min-max 3}, we first observe that $u\wedge v\in\C^n$ implies 
$$
m(u\vee v)=m(u)\vee m(v)
$$
(see the next lemma). Let $t>m(u\vee v)$; then, by Corollary \ref{corollario sottolivelli} and \eqref{min-max 1}:
\begin{eqnarray*}
\cl(\{u< t\})\cap\cl(\{v< t\})=\{u\le t\}\cap\{v\le t\}=\{u\vee v\le t\}=\cl(\{u\vee v\le t\}).
\end{eqnarray*}
If we assume that $t\le m(u\vee v)$, then we have $t\le m(u)$ or $t\le m(v)$. In the first case $\{u<t\}=\emptyset$, so that the left hand-side
of \eqref{min-max 3} is empty. On the other hand $u\vee v\le u$ implies that the right hand-side is empty as well. The case $t\le m(v)$ is completely 
analogous.
\end{proof}

\begin{lemma} If $u,v\in\C^n$ are such that $u\wedge v\in\C^n$, then
$$
m(u\vee v)=m(u)\vee m(v).
$$
\end{lemma}

\begin{proof} The inequality $m(u\vee v) \ge m(u)\vee m(v)$ is obvious. To prove the reverse inequality, let $t\ge m(u)\vee m(v)$; hence
$$
\{u\le t\}\ne\emptyset,\quad \{v\le t\}\ne\emptyset
$$
and 
$$
\{u\le t\}\cup\{v\le t\}=\{u\wedge v\le t\}\in\K^n,
$$
where the last relation comes from the assumption $u\wedge v\in\C^n$. Hence $\{u\le t\}$ and $\{v\le t\}$ are non-empty convex bodies such that their union is also
a convex body. This implies that they must have a non-empty intersection. But then
$$
\{u\vee v\le t\}=\{u\le t\}\cap\{v\le t\}\ne\emptyset,
$$
i.e. $m(u\vee v)\le t$.
\end{proof}

\bigskip

We conclude this section with a proposition (see \cite[Lemma 2.5]{Colesanti-Fragala}) which will be frequently used throughout the paper.

\begin{proposition}\label{limitazione-uniforme}
If $u\in\C^n$ there exist two real numbers $a$ and $b$, with
$a>0$, such that
$$
\mbox{$u(x)\ge a|x|+b$ for every $x\in\R^n$.}
$$
\end{proposition}

\section{Valuations on $\C^n$}

\begin{definition} A valuation on $\C^n$ is a map $\mu\,:\,\C^n\to\R$ such that
$\mu(\boldsymbol{\infty})=0$ and 
$$
\mu(u\vee v)+\mu(u\wedge v)=
\mu(u)+\mu(v)
$$
for every $u,v\in\C^n$ such that $u\wedge v\in\C^n$. 

A valuation $\mu$ is said:
\begin{itemize}
\item {\bf\em rigid motion invariant}, if $\mu(u)=\mu(u\circ T)$ for every $u\in\C^n$ and for every 
rigid motion $T$;
\item {\bf\em monotone decreasing} (or just monotone), if $\mu(u)\le\mu(v)$ whenever $u,v\in\C^n$ and $u\ge v$ point-wise in 
$\R^n$;
\item {\bf\em $k$-simple} ($k\in\{1,\dots,n\}$) if $\mu(u)=0$ for every $u\in\C^n$ such that $\dim(\dom(u))<k$;  
\item{\bf\em simple}, if $\mu$ is $n$-simple, i.e. if $\mu(u)=0$ for every $u$ such that $\dom(u)$ has no interior points.
\end{itemize}
\end{definition}

The following simple observation will turn out to be very important.

\label{aliter}
\begin{remark}
{\em Every monotone decreasing valuation $\mu$ on $\C^n$ is non-negative.
If we set $\boldsymbol{\infty}(x)=\infty$ for all $x\in\R^n$, then
$u\le\boldsymbol{\infty}$ holds for each $u\in\C^n$, which in turn leads to $\mu(u)\ge\mu(\boldsymbol{\infty})=0$ by monotonicity.} 
\end{remark}

In the sequel other features of valuations will be considered, like {\em monotone-continuity} and homogeneity. Concerning the 
latter, the definition is the following. 

\begin{definition}\label{hom-val}
Let $\mu$ be  valuation on $\C^n$ and let $\alpha\in\R$; we say that $\mu$ is positively homogeneous of order $\alpha$, or simply $\alpha$-homogeneous, 
if for every $u\in\C^n$ and every $\lambda>0$ we have
$$
\mu(u)=\lambda^\alpha \mu(u_\lambda)
$$
where $u_\lambda\,:\,\R^n\to\R\cup\{\infty\}$ is defined by
$$
u_\lambda(x)=u\left(\frac x\lambda\right)\quad\forall\, x\in\R^n
$$
(note that $u\in\C^n$ implies $u_\lambda\in\C^n$).
\end{definition}

\begin{remark} {\rm Other definitions of homogeneous valuations are possible. For instance one could consider valuations for which there 
exists $\alpha\in\R$ such that
$$
\mu(\lambda u)=\lambda^\alpha\mu(u),\quad\forall\, u\in\C^n,\,\forall\,\lambda>0.
$$ 
This corresponds to homogeneity with respect to a vertical stretching of the graph of $u$, while Definition \ref{hom-val} involves a horizontal
stretching. In addition, one could consider a more general type of homogeneity where both types of dilations (vertical and horizontal) are 
simultaneously in action. Definition \ref{hom-val} is more natural from the point of view of convex bodies. Indeed, if $u=I_K$ with $K\in\K^n$, 
then $u_\lambda$ is the indicatrix function of the dilated body $\lambda K$.}
\end{remark}

The next one is the definition of monotone-continuous valuations.

\begin{definition}\label{m-continuous} Let $\mu$ be a valuation on $\C^n$; $\mu$ is called monotone-continuous, or simply
m-continuous, if the following property is verified: given a sequence  $u_i\in\C^n$, $i\in\N$, and $u\in\C^n$, such that: 
$$
u_i\ge u_{i+1}\ge u \quad\mbox{in $\R^n$, for every $i\in\N$}, 
$$
and
$$
\lim_{i\to\infty} u_i(x)=u(x)\quad\forall\, x\in{\rm relint}(\dom(u))
$$
we have
$$
\lim_{i\to\infty}\mu(u_i)=\mu(u).
$$
\end{definition}

We recall that ``relint'' denotes the relative interior of a convex set (see the definition given in section \ref{subsection convex bodies}).

\begin{remark}\label{continuita' su successioni monotone di cc}
{\em Let $K_i$, $i\in\N$, be a sequence converging to a convex body $K$ in the Hausdorff metric. Assume moreover that
the sequence is monotone increasing:
$$
K_i\sub K_{i+1}\sub K\quad\forall\,i\in\N.
$$
Then the corresponding sequence of indicatrix functions $I_{K_i}$, $i\in\N$, is decreasing, it verifies $I_{K_i}\ge I_K$ point-wise in
$\R^n$, for every $i$, and
$$
\lim_{i\to\infty} I_{K_i}(x)=I_K(x)\quad\forall\,x\in{\rm relint}(K).
$$
Hence, if $\mu$ is an m-continuous valuation on $\C^n$,
$$
\lim_{i\to\infty}\mu(I_{K_i})=\mu(I_K).
$$}
\end{remark}




\section{Geometric densities}\label{sec geometric densities} 

Throughout this section, $\mu$ will be a rigid motion invariant and monotone decreasing valuation on $\C^n$. 

Let $t\in\R$ be fixed, and consider the following application $\sigma_t$ defined on $\K^n$:
$$
\sigma_t\,:\,\K^n\to\R,\quad
\sigma_t(K)=\mu(t+I_K)\quad\forall\, K\in\K^n.
$$
It is straightforward to check that $\sigma_t$ is a rigid motion invariant valuation on $\K^n$. Moreover, if $K\sub L$, then 
$I_K\ge I_L$, so that $\sigma_t(K)\le\sigma_t(L)$, i.e. $\sigma_t$ is monotone increasing.  By Theorem \ref{Hadwiger-bodies} there 
exist $(n+1)$ non-negative coefficients, that we will denote by $f_0,f_1,\dots, f_n$, such that
$$
\sigma_t(K)=\sum_{k=0}^n f_k\,V_k(K)\quad\forall\, K\in\K^n.
$$
The numbers $f_k$'s clearly depend on $t$, i.e. they are real-valued functions defined on $\R$; we will refer to these functions as the
{\em geometric densities} of $\mu$. 

We prove that the monotonicity of $\mu$ implies that these functions are monotone decreasing.
Fix $k\in\{0,\dots,n\}$ and let $B^k$ be a closed $k$-dimensional ball of radius $1$ in $\R^n$; note that
$$
V_j(B^k)=0\quad\forall\, j=k+1,\dots,n,
$$
while
$$
V_{k}(B^k)=:\kappa_k=\mbox{(the $k$-dimensional volume of $B^k$)}>0.
$$
Fix $r\ge0$; $V_j$ is positively homogeneous of order $j$, hence 
for every $t\in\R$ we have
$$
\mu(t+I_{rB^k})=\sum_{j=1}^k r^{j}V_j(B^k)f_j(t).
$$
Hence we get
$$
f_k(t)=V_k(B^k)\cdot\lim_{r\to \infty}\frac{\mu(t+I_{rB^k})}{r^{k}}.
$$
On the other hand, as $\mu$ is decreasing, the function $t\to\mu(t+I_{rB^k})$ is decreasing for every $r\ge0$; this proves that
$f_k$ is decreasing.

\begin{proposition}\label{geometric-densities} Let $\mu$ be a rigid motion invariant and decreasing valuation defined on $\C^n$. Then there
exists $(n+1)$ functions $f_0,f_1,\dots,f_n$, defined on $\R$, non-negative and decreasing, such that for every convex body $K\in\K^n$
and for every $t\in\R$
\begin{equation}\label{eq. geom. dens.}
\mu(t+I_K)=\sum_{k=0}^nf_k(t)V_k(K).
\end{equation}
\end{proposition}

%
%

If in addition to the previous assumption the valuation $\mu$ is m-continuous, then all its geometric densities are {\em right-continuous},
i.e.
$$
\lim_{t\to t_0^+}f_i(t)=f_i(t_0)\quad\forall t_0\in\R,\,i\in\{0,\dots,n\}.
$$
Indeed, for every convex body $K$ the function 
$$
t\mapsto \mu(t+I_K)
$$
is right-continuous, by the definition of m-continuity. If we chose $K=\{0\}$, as $V_k(K)=0$ for $k\ge1$ and $V_0(K)=1$,
we have, by \eqref{eq. geom. dens.}
$$
\mu(t+I_K)=f_0(t)\quad\forall\,t\in\R.
$$
This proves that $f_0$ is right-continuous. If we now take $K$ to be a one-dimensional convex body, such that $V_1(K)=1$, we have 
that $V_k(K)=0$ for every $k\ge2$, hence
$$
\mu(t+I_K)=f_0(t)+f_1(t)\quad\forall\,t\in\R.
$$
As the left hand-side is right-continuous and $f_0$ is also right-continuous (by the previous step) then $f_1$ must have the same property. Proceeding 
in a similar way we obtain that each $f_k$ is right-continuous.

\begin{proposition}\label{continuity of geometric densities} 
Let $\mu$ be a rigid motion invariant, monotone and m-continuous valuation on $\C^n$. Then its geometric densities
$f_i$, $i\in\{1,\dots,n\}$, are right-continuous in $\R$.
\end{proposition}

Assume now that $\mu$ is positively homogeneous of some order $\alpha$; then it is readily checked that for every $t\in\R$ the valuation 
$\sigma_t$ defined at the beginning of this section is positively homogeneous of the same order, i.e.
$$
\sigma_t(s\cdot K)=s^\alpha\sigma_t(K)\quad\forall\, K\in\K^n,\, s\ge0.
$$
On the other hand, each $\sigma_t$ is a linear combination of the intrinsic volumes $V_k$'s, and  $V_k$ is positively homogenous of order $k$. 
We are led to the following conclusion.

\begin{corollary}\label{order of hom.} Let $\mu$ be a rigid motion invariant and monotone decreasing valuation on $\C^n$ and assume that it is 
$\alpha$-homogeneous for some $\alpha\in\R$. Then necessarily $\alpha\in\{0,1,\dots,n\}$ and $f_k\equiv0$ for every $k\ne\alpha$.
\end{corollary}

We are in position to prove a characterization result for $0$-homogeneous valuations which are also monotone and m-continuous. We recall that, 
for $u\in\C^n$, $m(u)$ is the minimum of $u$ on $\R^n$.

\begin{proposition}\label{characterization of zero homogeneous valuations}   Let $\mu$ be a rigid motion invariant, monotone decreasing
and continuous valuation on $\C^n$ and assume that it is $0$-homogeneous. Then, for every $u\in\C^n$ we have
$$
\mu(u)=f_0(m(u)). 
$$
\end{proposition}

\begin{proof} We first prove the claim of this proposition under the additional assumption that $\dom(u)$ is bounded; let $K$ be a convex body containing
$\dom(u)$. Moreover, let $x_0\in\R^n$ be such that $u(x_0)=m(u)$. Then
$$
m(u)+I_K(x)\le u(x)\le
m(u)+I_{\{x_0\}}(x)\quad\mbox{for every $x$ in $\R^n$.}
$$
As $\mu$ is monotone decreasing
$$
\mu(m(u)+I_K)\ge \mu(u)\ge
\mu(m(u)+I_{\{x_0\}}).
$$
On the other hand, by Poposition \ref{geometric-densities} and Corollary \ref{order of hom.} we have
$$
\mu(m(u)+I_{\{x_0\}})=\mu(m(u)+I_K)=f_0(m(u)),
$$
hence $\mu(u)=f_0(m(u))$. To extend the result to the general case, for $u\in\C^n$ and $i\in\N$ let
$$
u_i=u+I_{B(i)}
$$
where $B(i)$ denotes the closed ball of radius $i$ centred at the origin. The sequence $u_i$ is contained in $\C^n$, is monotone decreasing and 
converges point-wise to $u$ in $\R^n$. As $\mu$ is m-continuous we have, by the previous part of the proof,
$$
\mu(u)=\lim_{i\to\infty}\mu(u_i)=\lim_{i\to\infty}f_0(m(u_i)).
$$
On the other hand, as $m(u)=\min_{\R^n}u$, by the point wise convergence we have that for $i$ sufficiently large $m(u_i)=m(u)$.
\end{proof}

Another special case in which more information can be deduced on geometric densities, is when the valuation $\mu$ is simple.

\begin{proposition}\label{geometric densities simple valuations}
Let $\mu$ be a rigid motion invariant, monotone decreasing and simple valuation on $\C^n$. Then, for each $k\in\{0,\dots,n-1\}$,
the $k$-th geometric density $f_k$ of $\mu$ is identically zero.
\end{proposition}

\begin{proof} 
Fix $t\in\R$; the valuation $\sigma_t :\mathcal{K}^n \to\R$ defined by 
$$ \sigma_{K}=\mu(t+I_K)$$ 
is monotone, rigid motion invariant and simple; the volume theorem (Corollary \ref{volt}) and the definition of geometric densities yield
$$\sigma_t= f_n(t) V_n.$$
In other words, $f_k(t)=0$ for every $k=0,\dots,n-1$ and $t\in\R$.
\end{proof}

The following result relates homogeneity and simplicity, and its proof makes use of geometric densities.

\begin{proposition} \label{k-homogeneous implies k-simple} Let $\mu\,:\,\C^n\to\R$ be a valuation with the following properties:
\begin{itemize}
\item $\mu$ is rigid motion invariant;
\item $\mu$ is monotone decreasing;
\item $\mu$ is $k$-homogeneous;
\item $\mu$ is m-continuous.
\end{itemize}
Then $\mu$ is $k$-simple.
\end{proposition}
\begin{proof} Let $f_0,f_1,\dots,f_n$ be the geometric densities of $\mu$. As $\mu$ is $k$-homogeneous,
$f_i\equiv0$ for every $i\ne k$. Let $u\in\C^n$ be such that $\dim(\dom(u))<k$. For $i\in\N$ let
$Q_i=[-i,i]^n=[-i,i]\times\dots\times[-i,i]$, and set
$$
u_i=u+I_{Q_i}.
$$
Clearly $\dom(u_i)=\dom(u)\cap Q_i$; in particular $\dim(\dom(u_i))<k$ for every $i$. Let
$$
m_i=m(u_i),\quad\,i\in\N.
$$
As $\mu$ is monotone we have
$$
0\le\mu(u_i)\le\mu(m_i+I_{\cl(\dom(u_i))})=f_k(m_i) V_k(\cl(\dom(u_i))).
$$
On the other hand $V_k(\cl(\dom(u_i)))=0$, as $\dim(\dom(u_i))<k$. Hence $\mu(u_i)=0$ for every $i$.
To conclude, note that $u_i$, $i\in\N$, is a decreasing sequence of functions in $C^n$, converging to $u$ point-wise
in $\R^n$. As a consequence of m-continuity of $\mu$ we have
$$
\mu(u)=\lim_{i\to\infty}\mu(u_i)=0.
$$
\end{proof}

\subsection{Regularization of geometric densities}\label{regularization of geometric densities}

In the sequel sometimes it will be convenient to work with valuation having geometric densities with more regularity than that of a decreasing function. In this section 
we describe a procedure which allows to approximate (in a suitable sense) a valuation with a sequence of valuations having smooth densities.

Let $g\,:\,\R\to\R$ be a standard mollifying kernel, i.e. $g$ has the following properties:
$g\in C^\infty(\R)$, $g(t)\ge0$ for every $t\in\R$, the support of $g$ is contained in $[-1,1]$ and
$$
\int_{\R} g(t)dt=1.
$$ 
For $\epsilon>0$ let $g_\epsilon\,:\,\R\to\R$ be defined by
$$
g_\epsilon(t)=\frac1\epsilon \,g\left(
\frac t\epsilon
\right).
$$
Then $g_\epsilon\in C^\infty(\R)$; $g_\epsilon(t)\ge0$ for every $t\in\R$; the support of $g_\epsilon$ is contained in $[-\epsilon,\epsilon]$ and 
$$
\int_{\R} g_\epsilon(t)dt=1.
$$ 

Now let $u\in\C^n$ and consider the function $\phi\,:\,\R\to\R$ defined by
$$
\phi(t)=\mu(u+t).
$$
By the properties of $\mu$, this is a non-negative and decreasing function. For $\epsilon>0$ and $t\in\R$ set:
$$
\phi_\epsilon(t)=(\phi\star g_\epsilon)(t)=\int_{\R}\phi(t-s) g_\epsilon(s)ds
=\int_{\R}\mu(u+t-s) g_\epsilon(s)ds.
$$
Then $\phi_\epsilon$ is a non-negative decreasing function of class $C^\infty(\R)$. Moreover, by the properties of the kernel $g$,
$$
\lim_{\epsilon\to0^+}\phi_\epsilon(t)=\phi(t)
$$ 
for every $t\in\R$ where $\phi$ is continuous (in fact, for every Lebesgue point $t$ of $\phi$, see for instance \cite{EG}); 
in particular $\phi_{\epsilon}\to\phi$ a.e. in $\R$ as $\epsilon\to0^+$.

For $\epsilon>0$ we define $\mu_\epsilon\,:\,\C^n\to\R$ as
$$
\mu_\epsilon(u)=\phi_\epsilon(0)=\int_{\R}\phi(-s)g_\epsilon(s)ds=
\int_{\R}\mu(u-s)g_\epsilon(s)ds.
$$
It is a straightforward exercise to verify that $\mu_\epsilon$ inherits most of the properties of $\mu$: it is a valuation,
rigid motion invariant, non-negative and decreasing. Moreover, if $f_{k,\epsilon}$, for $k\in\{0,\dots,n\}$, denote the 
geometric densities of $\mu_\epsilon$, we have
$$
f_{k,\epsilon}=f_k\star g_\epsilon
$$
for every $k\in\{0,\dots,n\}$, where $f_0\dots,f_n$ are the densities of $\mu$. Indeed, for every convex body $K\in\K^n$ and every $t\in\R$ 
we have
\begin{eqnarray*}
\mu_\epsilon(t+I_K)&=&\int_{\R}\mu(t-s+I_K)g_\epsilon(s)ds
=\int_{\R}\sum_{k=0}^nf_k(t-s)V_k(K)g_\epsilon(s)ds\\
&=&\sum_{k=0}^n(f_k\star g_\epsilon)(t)V_k(K).
\end{eqnarray*}

The core properties of the above construction can be summed up in the following proposition.

\begin{proposition}\label{mollification}
Let $\mu\,:\,\C^n\to\R$ be a rigid motion invariant and decreasing valuation. Then there exists a family of 
rigid motion invariant and decreasing valuations $\mu_\epsilon$, $\epsilon>0$, such that:
\begin{itemize}
\item[1.] for every $\epsilon>0$ the geometric densities $f_{0,\epsilon},\dots,f_{n,\epsilon}$ of $\mu_\epsilon$ belong to $C^\infty(\R)$;
\item[2.] for every $k\in\{0,\dots,n\}$, $f_{k,\epsilon}\to f_k$ a.e. in $\R$, as $\epsilon\to0^+$, where $f_0,\dots,f_n$ are the geometric densities of $\mu$;
\item[3.] for every $u\in\C^n$, $\mu_\epsilon(u+t)\to\mu(u+t)$ for a.e. $t\in\R$, as $\epsilon\to0^+$.
\end{itemize}
\end{proposition}

\section{Integral valuations}\label{sec integral valuations}
In this section we introduce a class of integral valuations which will turn out to be crucial in the characterization results that we will present in the sequel. 
As we will see, they are similar to those introduced by Wright in \cite{Wright} and subsequently studied by Baryshnikov, Ghrist and Wright in \cite{BGW}. 

Let $\nu$ be a (non-negative) Radon measure on the real line $\R$ and fix $k\in\{0,\dots,n\}$. For every $u\in\C^n$ we set
\begin{equation}\label{Wright-val}
\mu(u):=\int_\R V_k(\cl(\Omega_t))d\nu(t),
\end{equation}
where $\Omega_t=\{u<t\}$ for every $t\in\R$. As noted in sub-section \ref{On the intrinsic volumes of sub-level sets}, 
the function $t\mapsto V_k(\cl(\Omega_t))$ vanishes on $(-\infty,m(u)]$ and, for $k\ne0$ its $k$-th root
is concave in $(m(u),\infty)$, while for $k=0$ it is simply constantly 1 in $(m(u),\infty)$; hence it is Borel measurable. Moreover it is 
non-negative, so that it is integrable with respect to $\nu$. On the other hand its integral \eqref{Wright-val} might be $\infty$.  We first find equivalent
conditions on $\nu$ such that \eqref{Wright-val} is finite for every $u\in\C^n$.

\begin{proposition}\label{integrability condition} 
Let $\nu$ be a non-negative Radon measure on the real line. The integral \eqref{Wright-val} is finite for every
$u\in\C^n$ if and only if 
\begin{equation}\label{int-cond}
\int_{(0,\infty)} t^{k}d\nu(t)<\infty.
\end{equation}
\end{proposition}

\begin{proof} Assume that $\mu(u)$ is finite for every $u$. Choosing in particular $u\in\C^n$ defined by $u(x)=|x|$ we have that $V_k(\cl(\Omega_t))$ is zero for every
$t\le 0$. For $t>0$, $\cl(\Omega_t)$ is a ball centred at the origin with radius $t$, hence $V_k(\cl(\Omega_t))=c(n,k) t^k$ with $c(n,k)>0$. Therefore
$$
\mu(u)=\int_{(0,\infty)} c(n,k) t^k d\nu(t)<\infty.
$$
Vice versa, assume that condition \eqref{int-cond} holds. Given $u\in\C^n$ there exists $a>0$ and $b\in\R$ such that $u(x)\ge a|x|+b$ for every
$x$ (see Proposition \ref{limitazione-uniforme}). Hence, for $t\in\R$, $t\ge b$,
$$
\Omega_t\sub\{x\in\R^n \,:\, a|x|+b<t\}=\left\{x\in\R^n \,:\,|x|<\frac{t-b}{a} \right\},
$$
while $\Omega_t$ is empty for $t\le b$. By the monotonicity of intrinsic volumes
$$
\mu(u)=\int_\R V_i(\cl(\Omega_t))d\nu(t)\le
c(n,k) \int_{(b,\infty)}\left(
\frac{t-b}{a}
\right)^k d\nu(t),
$$
and the last integral is finite by \eqref{int-cond}.
\end{proof}

From now on in this section the summability condition \eqref{int-cond} will always be assumed (with respect to some fixed index $k\in\{0,\dots,n\}$).

\begin{proposition}\label{atomo} Let $k\in\{0,1,\dots,n\}$ and let $t$ be a fixed real number. Then the application
$u\mapsto V_k(\cl(\{u<t\}))$ 
\begin{itemize}
\item[{\it i)}] is a valuation;
\item[{\it ii)}] is rigid motion invariant;
\item[{\it iii)}] is monotone;
\item[{\it iv)}] is $k$-homogeneous;
\item[{\it v)}] is m-continuous.
\end{itemize}
\end{proposition}
\begin{proof}
{\em i)} 
The condition on $\boldsymbol{\infty}$ is easily verified, as a matter of fact
$$
V_k(\cl(\{\boldsymbol{\infty}<t\}))=V_k(\emptyset)=0. 
$$
Let now $u,v\in\C^n$ be such that $u\wedge v\in\C^n$. By Proposition \ref{proposition min-max}, 
for every $t\in\R$ we have
$$
\cl(\{u\wedge v< t\})=\cl(\{u< t\})\cup\cl(\{v< t\}),\quad
\cl(\{u\vee v< t\})=\cl(\{u< t\})\cap\cl(\{v< t\}).
$$
Consequently, as intrinsic volumes are valuations, we get
\begin{eqnarray*}
V_k(\cl(\{u\wedge v< t\}))&+&V_k(\cl(\{u\vee v< t\}))=\\
&=&V_k(\cl(\{u< t\})\cup\cl(\{v< t\}))+V_k(\cl(\{u< t\})\cap\cl(\{v< t\}))\\
&=&V_k(\cl(\{u< t\}))+V_k(\cl(\{v< t\})).
\end{eqnarray*}

{\em ii)} Let $u\in\C^n$ and let $T\,:\,\R^n\to\R^n$ be a rigid motion; let moreover $v:=u\circ T\in\C^n$, i.e.
$v(y)=u(T(y))$ for every $y\in\R^n$.  Then, for every $t\in\R$,
$$
\{y\,:\,v(y)<t\}=\{y\,:\,u(T(y))< t\}=
T^{-1}(\{x\,:\,u(x)< t\}).
$$
As intrinsic volumes are invariant with respect to rigid motions, we have
$$
V_k(\{u<t\})=V_k(\{v<t\}),
$$
for every $t$.

{\em iii)} As for monotonicity,
if $u,v\in\C^n$ and $u\le v$ point-wise on $\R^n$, then for every $t\in\R$,
$$
\{v< t\}\sub\{u< t\},
$$
and thus 
$$
\cl (\{v<t\}) \sub \cl (\{u<t\}).
$$

By the monotonicity of intrinsic volumes
$$
V_k(\cl(\{v< t\}))\le V_k(\cl(\{u<t\})).
$$

{\em iv)} Let $u\in\C^n$ and $\lambda>0$, and define $u_\lambda$ by
$$
u_\lambda(x)=u\left(\frac x\lambda\right).
$$
For  $t\in\R$ we have
$$
\{x\,:\,u_\lambda(x)< t\}=\left\{
x\,:\,u\left(\frac x\lambda\right)< t
\right\}=
\lambda\{x\,:\,u(x)< t\}.
$$
Then, by homogeneity of intrinsic volumes
$$
V_k(\cl(\{u_\lambda<t\}))=V_k(\lambda\,\cl(\{u< t\}))
=\lambda^kV_k(\cl(\{u< t\})).
$$
{\em v)} Let $u\in\C^n$ and let $u_i$, $i\in\N$, be a sequence in $\C^n$, point-wise decreasing
and converging to $u$ in the relative interior of $\dom(u)$. 
We want to prove that 
$$\lim_{i\to\infty} V_k(\cl(\{u_i<t\}))=V_k(\cl(\{u<t\})).$$
Let $j\in\{0,1,\dots, n\}$ be the dimension of $\dom(u)$; for every $i\in\N$, as $u_i\ge u$ we have that
$\dom(u_i)\sub\dom(u)$, and, in particular, the dimension of the domain of $u_i$ is less than or equal to $j$. If $k>j$, then 
$V_k(\cl(\{u< t\}))=0$ for every $t\in\R$ and analogously $V_k(\cl(\{u_i< t\}))=0$ for every $i$, so that the assert of the proposition
holds true. Hence we may assume that $k\le j$ and, up to restricting all involved functions to a $j$-dimensional affine subspace of $\R^n$
containing the domain of $u$, we may assume without loss of generality that $j=n$. 

As usual, we denote $\min_{\R^n}u$ by $m(u)$. 
If $t\le m(u)$, then, for all $i$, $\{u_i<t\}=\{u<t\}=\emptyset$ and the claim holds trivially.
Let now $t>m(u)$, then, by Corollary \ref{corollario sottolivelli}
$$
\cl(\{u<t\})=K_t:=\{u\le t\}.
$$
As $\dim(\dom(u))=n$ and $t>m(u)$, $K_t$ is a convex body with non-empty interior (this follows from the convexity of $u$). Let
$$
\Omega^i_t=\{u_i<t\},\quad 
K^i_t=\cl(\Omega^i_t),\quad
i\in\N.
$$
Clearly $K^i_t\sub K_t$ for every $i$. On the other hand, if $x$ is an interior point of $K_t$, then $u(x)<t$ (see Lemma \ref{lemma sottolivelli}). 
Hence $u_i(x)<t$ for sufficiently large $i$, which leads to
$$
\bigcup_{i\in\N} K^i_t\supseteq\interno(K_t).
$$
As a consequence, $K^i_t$ converges to $K_t$ in the Hausdorff metric, and then (by the continuity of intrinsic volumes)
$$
\lim_{i\to\infty}V_k(K^i_t)=V_k(K_t), \quad\forall\, t>m(u),
$$
as we wanted.
\end{proof}

\begin{corollary}\label{m-continuita' e continuita' destra} Let $k\in\{0,1,\dots,n\}$ and let $\nu$ be a Radon measure on $\R$ which verifies condition \eqref{int-cond}. Then the application
$\mu\,:\,\C^n\to\R$ defined by
$$
\mu(u)=\int_{\R}V_k(\cl(\{u<t\}))d\nu(t),\quad u\in\C^n
$$
has the following properties:

\begin{itemize}
\item[{\it i)}] it is a valuation;
\item[{\it ii)}] it is rigid motion invariant;
\item[{\it iii)}] it is monotone;
\item[{\it iv)}] it is $k$-homogeneous;
\item[{\it v)}] it is m-continuous.
\end{itemize}
\end{corollary}
 
\begin{proof}
Claims {\it i)} - {\it iv)} follow easily from Proposition \ref{atomo} by integration.
The proof of the m-continuity of $\mu$ is a bit more delicate.
Let $u\in\C^n$ and let $u_i$, $i\in\N$, be a sequence in $\C^n$, point-wise decreasing
and converging to $u$ in $\relint(\dom(u))$.
 

As $u_i\ge u$ we have, for every $t\in\R$, 
$$
\{u_i < t\}\sub\{u< t\}\,\Rightarrow\,
V_k(\cl(\{u_i< t\}))\le V_k(\cl(\{u< t\}))
$$
for every $i\in\N$.

By Proposition \ref{atomo} we know that
$$
\lim_{i\to\infty}V_k(\cl(\{u_i<t\}))=V_k(\cl(\{u<t\})), \quad\forall\, t.
$$
This fact and the monotone convergence theorem imply
$$
\lim_{i\to\infty}\int_{\R} V_k(\cl(\{u_i< t\}))d\nu(t)
=\int_{\R} V_k(\cl(\{u<t\}))d\nu(t).
$$
\end{proof}

Let $\mu$ be a valuation of the form \eqref{Wright-val}; by Proposition \ref{geometric-densities} and Corollary \ref{order of hom.}, $\mu$ has exactly
one geometric density which is not identically zero, i.e. $f_k$; this can be explicitly computed in terms of the measure $\nu$. Let $K\in\K^n$ be such that
$V_k(K)>0$, then, for $t\in\R$
\begin{equation*}
\begin{split}
\mu(t+I_K)&=
f_k(t)V_k(K)=\int_\R V_k(\cl (\{x\,:\,t+I_K(x)<s\}))d\nu(s) \\ 
&=V_k(K)\int_{\{s \,:\, t< s\}}d\nu(s)=V_k(K)\int_{(t,\infty)} d\nu(s);
\end{split}
\end{equation*}
i.e.
\begin{equation}\label{geometric density of an integral valuation}
f_k(t)
=\nu((t,\infty))
\quad\forall\,t\in\R.
\end{equation}

We observe that this is a non-increasing function and, by the basic properties of measures, it is right-continuous.

\subsection{An equivalent representation formula}

As in the previous part of this section, $\nu$ will be a non-negative Radon measure on $\R$ verifying the integrability condition
\eqref{int-cond}, where $k$ is a fixed integer in $\{0,1,\dots,n\}$. Moreover, $f\,:\,\R\to\R$ is defined by
\begin{equation}\label{definition of f}
f(t)=\int_{(t,\infty)} d\nu(t)=\nu((t,\infty)).
\end{equation}
We first consider the case $k\ge 1$. Note that \eqref{int-cond} is equivalent to
\begin{equation}\label{int-cond-f}
\int_0^\infty t^{k-1}f(t)dt<\infty.
\end{equation}
Indeed
\begin{eqnarray*}
\int_0^\infty t^{k-1}f(t)dt&=&\int_0^\infty t^{k-1}\int_{(t,\infty)} d\nu(s)dt
=\int_{(0,\infty)}\int_0^s t^{k-1}dt d\nu(s)
=\frac 1k\int_{(0,\infty)} s^k d\nu(s). 
\end{eqnarray*}

Now let $u\in\C^n$ and let $v_k(u;\cdot)$ be the function defined in section \ref{On the intrinsic volumes of sub-level sets}
$$
v_k(u;t)=V_k(\cl(\{u< t\}))\quad t\in\R.
$$
For simplicity we set $h(t)=v_k(u;t)$ for every $t\in\R$; 
$h$ is a monotone non-decreasing function identically vanishing on $(-\infty,m(u)]$ and $h^{1/k}$ is concave in $(m(u),\infty)$, as pointed out in 
section \ref{On the intrinsic volumes of sub-level sets}; in particular $h$ is locally Lipschitz. This implies 
(see for instance \cite[Chapter 3]{AFP}) that the product $fh$ is a function of bounded
variation in $(m(u),\infty)$ and its weak derivative is the measure
$$
-h\nu+h'f{\mathcal H}^1
$$
(we recall that ${\mathcal H}^1$ is the one-dimensional Lebesgue measure). Note also that $fh$ is right-continuous, as
$f$ has this property. Hence for every $t_0,t\in\R$, with $m(u)< t_0\le t$,
$$
f(t)h(t)=f(t_0)h(t_0)+\int_{(t_0,t)}h'(s)f(s)ds-\int_{(t_0,t)}h(s)d\nu(s).
$$
If we let $t_0\to m(u)^+$ we get
\begin{equation}\label{proof2}
f(m(u))V_k(\{u=m(u)\})+\int_{(m(u),t)}f(s)h'(s)ds=
f(t)h(t)+\int_{(m(u),t)} f(s)d\nu(t).
\end{equation}
Indeed, as we proved in section \ref{On the intrinsic volumes of sub-level sets},
$$
\lim_{t\to m(u)^+} h(t)=V_k(\{u= m(u)\}).
$$

By Lemma \ref{limitazione-uniforme}, there is a constants $a>0$ such that
$$
h(t)\le at^i
$$
for $t$ sufficiently large. Hence $h(t) f(t)\le a t^i f(t)$. On the other hand, the integrability condition \eqref{int-cond-f} implies
$$
\liminf_{t\to\infty} t^i f(t)=0,
$$
so that
$$
\liminf_{t\to\infty}f(t)h(t)=0.
$$
Hence, passing to the limit for $t\to\infty$ in \eqref{proof2}, we get
$$
f(m(u))V_k(\{u=m(u)\})+\int_{(m(u),\infty)} f(s)h'(s)ds=
\int_{(m(u),\infty)} h(s)d\nu(s)=
\int_{\R}h(s)d\nu(s)
$$
(the last equality is due to: $h\equiv0$ in $(-\infty,m(u)]$).
Recalling the structure of the weak derivative of the function $h$ proven in Proposition \ref{structure beta_i}, we may write
\begin{equation}\label{proof3}
\int_{\R}f(s)d\beta_k(u;s)=
\int_{\R}V_k(\cl(\{u< s\}))d\nu(s).
\end{equation}
The above formula is proven for $k\ge 1$. The case $k=0$ is straightforward, indeed
$$
V_0(\{u<t\})=
\left\{
\begin{array}{lll}
\mbox{$0$ if $t\le m(u)$,}\\
\\
\mbox{$1$ if  $t> m(u)$}
\end{array}
\right.
$$  
so that \eqref{proof3} becomes
$$
f(m(u))=\int_{(m(u),\infty)} d\nu(s),
$$
which is true by the definition of $f$. The previous considerations provide the proof of the following result.

\begin{proposition}\label{Wright representation} Let $k\in\{0,\dots,n\}$, let $\nu$ be a non-negative Radon measure on $\R$ verifying the integrability condition
\eqref{int-cond}, and let $f$ be defined as in \eqref{definition of f}. Let $\mu\,:\,\C^n\to\R$ be the valuation defined by \eqref{Wright-val}. Then 
for every $u\in\C^n$
$$
\mu(u)=\int_{\R}f(t)d\beta_k(u;t).
$$
\end{proposition}

Valuations expressed as in the above proposition were considered in \cite{Wright} and \cite{BGW}.

\bigskip

In the remaining part of this section we analyze two special cases of the integral valuation introduced so far, corresponding to the indices
$k=0$ and $k=n$, which can be written in a simpler alternative form.

\subsection{The case $k=0$} Let $\nu$ be a Radon measure on $\R$; the integrability condition \eqref{int-cond} for $k=0$ is just
$$
\int_0^\infty d\nu(s)<\infty,
$$
which is equivalent to saying that the function $f$ defined by \eqref{definition of f} is well defined (i.e. finite) in $\R$. Let $u\in\C^n$; as
we pointed out before,
$$
V_0(\cl(\{u<t\}))=
\left\{
\begin{array}{lll}
\mbox{$0$ if $t\le m(u)$,}\\
\\
\mbox{$1$ if $t> m(u)$.}
\end{array}
\right.
$$
Then
$$
\mu(u)=
\int_{(m(u),\infty)} d\nu(t)=f(m(u))
$$ 
for every $u\in\C^n$.

\subsection{The case $k=n$} 

\begin{proposition}\label{i=n} Let $\nu$ be a Radon measure on $\R$ verifying the integrability condition 
$$
\int_0^\infty t^nd\nu(t)<\infty,
$$
let $f$ be defined as in \eqref{definition of f}, and let $\mu$ be the valuation:
$$
\mu(u)=\int_\R V_n(\cl(\{u<t\}))d\nu(t),\quad u\in\C^n.
$$
Then
$$
\mu(u)=\int_{\dom(u)}f(u(x))dx\quad\forall\,u\in\C^n.
$$
\end{proposition}

\begin{proof} Let us extend $f$ to $\R\cup\{\infty\}$ setting $f(\infty)=0$.
As a direct consequence of the so-called layer cake (or Cavalieri's) principle and of the definition of $f$, we have that
$$
\int_{\R}{\mathcal H}^n(\{u<t\})d\nu(t)=
\int_{\R^n} f(u(x))dx=
\int_{\dom(u)}f(u(x))dx,
$$
where ${\mathcal H}^n$ denotes the Lebesgue measure in $\R^n$. 
On the other hand, for every $t\in\R$ the set $\{u<t\}$ is convex and bounded, so that its boundary is negligible with respect to the 
Lebesgue measure. Hence ${\mathcal H}^n(\{u<t\})=V_n(\cl(\{u<t\}))$ for every $t$.
\end{proof}

We can change the point of view and take the function $f$ as a starting point, instead of the measure $\nu$.

\begin{proposition}\label{Valutazioni integrali semplici}
Let $f\,:\,\R\to\R$ be non-negative, decreasing, and right-continuous. Define the mapping
$\mu\,:\,\C^n\to\R$ as follows
$$
\mu(u)=\int_{\dom(u)} f(u(x))dx
$$
for every $u\in\C^n$. Then:
\begin{itemize}
\item[{\it i})] $\mu$ is well defined (i.e. $\mu(u)\in\R$ for every $u\in\C^n$) if and only if 
$f$ verifies the following integrability condition
$$
\int_0^{\infty}t^{n-1}f(t)dt<\infty;
$$
\item[{\it ii})] $\mu$ is a valuation on $\C^n$, and it is rigid motion invariant, simple and decreasing;
\item[{\it iii})] $\mu$ is $n$-homogeneous;
\item[{\it iv})] $\mu$ is m-continuous.
\end{itemize}
\end{proposition}

The proof follows directly from the previous considerations and the dominated convergence theorem. A typical example in this sense is given by the application $\mu$ defined by
$$
\mu(u)=\int_{\dom(u)} e^{-u(x)}dx.
$$

\section{A decomposition result for simple valuations}\label{Partitions}

We start by introducing a particular class of partitions of convex bodies that will be used in the sequel.

\subsection{Convex partitions and inductive partitions}

\begin{definition} {\bf (Convex partition).} Let $K,K_1,\dots,K_N\in\K^n$; 
$$
{\mathcal P}:=\{K_1,\dots,K_N\}
$$
is called a convex partition of $K$ if
$$
K=\bigcup_{i=1}^N K_i\quad\mbox{and}\quad
\interno(K_i \cap K_j)=\emptyset\quad\forall\, i\ne j.
$$
\end{definition}

\begin{definition}{\bf (Inductive partition).} \label{indpart}
Let $K\in\K^n$. A convex partition $\mathcal{P}$ of $K$ is called an \emph{inductive partition} if there exists a sequence of convex bodies $H_1,\dots,H_l=K$ such that, for all $i=1,\dots,l$, one of the two following conditions holds true: 
\begin{itemize}
\item $H_i\in\mathcal{P}$,
\item $\exists j,k<i$ such that $H_i=H_j\cup H_k$ and $\interno(H_j\cap H_k)=\emptyset$.
\end{itemize}
\end{definition}

\bigskip

\noindent{\bf Examples.} 
It is immediate to prove that convex partitions made of one or two elements are inductive partitions.
On the other hand it is easy to construct a convex partition of three elements which is
not inductive. Let $K$ be a disk in the plane, centred at the origin, and consider three rays from the origin such that the angle between any
two of them is $2\pi/3$. These rays divide $K$ into three subsets $K_1$, $K_2$ and $K_3$ which form a convex partition $\mathcal P$ of $K$.
$\mathcal P$ is not an inductive partition.

\subsection{Complete partitions}\label{subsection Complete partitions}

From now on in this section $P$ will be a {\em polytope} of $\R^n$, i.e. the convex hull of finitely many points of $\R^n$. Note in particular
that $P\in\K^n$; moreover we will always assume that {\em $P$ has non-empty interior}.  We consider a convex partition 
$$
{\mathcal P}=\{P_1,\dots,P_N\}
$$
of $P$ whose elements are all polytopes, with non-empty interior; we will refer to such partitions as {\em polytopal} partitions.
 
If $H$ is a hyperplane (i.e. an affine subspace of dimension $(n-1)$) of $\R^n$, we can refine $\mathcal P$
by $H$ in the usual way. Let $H^+$ and $H^-$ be the closed half-spaces determined by $H$ and set
$$
P_i^+=P_i\cap H^+,\quad
P_i^-=P_i\cap H^-.
$$
Then 
$$
{\mathcal P}_H=\{P_i^+\,:\,i\in\{1,\dots,N\},\,\interno(P_i^+)\ne\emptyset\}\cup
\{P_i^-\,:\,i\in\{1,\dots,N\},\,\interno(P_i^-)\ne\emptyset\}
$$
is still a polytopal partition of $K$ that we indicate as the refinement of $\mathcal P$ by $H$.

\begin{definition} A polytopal partition is said to be complete if for every hyperplane $H$ which contains an $(n-1)$-dimensional 
facet of at least one element of $\mathcal P$ we have
$$
{\mathcal P}_H={\mathcal P}.
$$
\end{definition}

\begin{remark}\label{complete refinement}
{\em
Every polytopal partition $\mathcal P$ can be successively refined until it becomes, in a finite number of steps, a complete partition. Indeed, let
$\{H_1,\dots,H_R\}$ be the collection of all hyperplanes containing a $(n-1)$-dimensional facet of at least one element 
of $\mathcal P$. Then 
$$
(\dots(({\mathcal P})_{H_1})_{H_2}\dots)_{H_R}
$$
is a complete partition.}
\end{remark}


\begin{figure}[h]
    \centering
    \includegraphics[width=0.6\textwidth]{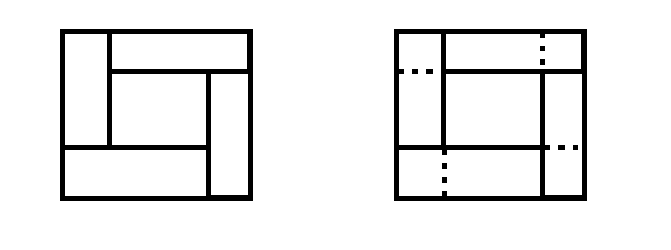}
    \caption{A partition being completed}
    \label{}
\end{figure}

\begin{proposition}\label{complete partitions are D-partitions}
Let $P$ be a polytope with non-empty interior and let ${\mathcal P}=\{P_1,\dots, P_N\}$ be a complete 
partition of $P$. Then $\mathcal P$ is an inductive partition of $P$.
\end{proposition}
\begin{proof}

The proof proceeds by induction on $N$. For $N=1$ the assert is true as any partition consisting of one element
is trivially an inductive partition. 
Let $N\ge 2$ and assume that the claim is true for every integer up to $(N-1)$. Let $\mathcal P$ be a complete
polytopal partition of $P$. Let $H$ be a hyperplane containing a $(n-1)$-dimensional facet of an element of
$\mathcal P$, intersecting the interior of $P$. Such a hyperplane exists because $P$ has non-empty interior and $N\ge 2$. 
Let $H^+$ and $H^-$ be the closed half-spaces determined by $H$. Then, as $\mathcal P$ is complete
(and each $P_i$ has non-emtpy interior), each $P_i$ is contained either in $H^+$ or in $H^-$.
Moreover, as $N\ge2$ and $H\cap\interno(P)\ne\emptyset$, there exist at least one element of $\mathcal P$ contained in $H^+$
and at least one element in $H^-$. Then clearly
\begin{eqnarray*}
{\mathcal P}^+&=&\mbox{$\{Q\in\mathcal{P} \,:\, Q\sub H^+\}$ is a complete partition of $P\cap H^+$,}\\
{\mathcal P}^-&=&\mbox{$\{Q\in\mathcal{P} \,:\, Q\sub H^-\}$ is a complete partition of $P\cap H^-$.}
\end{eqnarray*}
Each of these partitions has a number of elements which is strictly less than $N$. Consequently,
by the induction hypothesis, ${\mathcal P}^+$ is also an inductive partition of $P\cap H^+$ and 
${\mathcal H}^-$ is an inductive partition of $P\cap H^-$. 
Therefore, by Definition \ref{indpart}, there exist two sequences 
$$P_1^+,\dots, P_j^+=P\cap H^+ \quad \text{and} \quad P_1^-,\dots, P_k^-=P\cap H^-$$ that fulfill the required properties.  
We claim that such a sequence can be formed for the partition $\mathcal P$ as well: as a matter of fact consider the following 
$$ P_1^+,\dots,P_j^+,P_1^-,\dots,P_k^-,P.$$
As $P_j^+\cup P_k^-=P$ and $\interno(P_j^+\cap P_k^-)=\emptyset$ we conclude that $\mathcal P$ is an inductive partition too. 
\end{proof}

\subsubsection{Rectangular partitions}\label{Rectangular partitions}
A {\em rectangle} $R$ in $\R^n$ is a set of the form
$$
R=\{(x_1,\dots,x_n)\in\R^n\,:\,a_j\le x_j\le b_j\;\mbox{for every $j=1,\dots,n$}\},
$$
where, for $j=1,\dots,n$, $a_j$ and $b_j$ are real numbers such that $a_j< b_j$. In particular, $R$ is a convex polytope, and each of its facets
is parallel to a hyperplane of the form $e_j^\perp$, for some $j\in\{1,\dots,n\}$ (where $\{e_1,\dots,e_n\}$ is the canonical basis of $\R^n$). This property 
characterizes rectangles.

A {\em rectangular partition} of a rectangle $R$ is a partition
$$
{\mathcal P}=\{R_1,\dots,R_N\}
$$
of $R$ such that each $R_k$ is itself a rectangle. 

If ${\mathcal P}$ is a rectangular partition of a rectangle $R$, and we refine it so that its refinement ${\mathcal P}'$ is complete, as indicated in Remark 
\ref{complete refinement}, then ${\mathcal P}'$ is still a rectangular partition; indeed each facet of each element of ${\mathcal P}'$ is contained
in a hyperplane parallel to $e_j^\perp$, for some $j$.

\subsection{A decomposition result for simple valuations}

The following result is the main motivation for the definition of inductive partitions.

\begin{lemma}\label{decomposition of simple valuations} 
Let $\mu$ be a simple valuation on $\C^n$. Let $K$ be a convex body and let 
$$
{\mathcal P}=\{K_1,\dots,K_N\}
$$
be an inductive partition of $K$. Then, for every $u\in\C^n$,
$$
\mu(u+I_K)=
\sum_{i=1}^N 
\mu(u+I_{K_i}).
$$
\end{lemma}

\begin{proof} 
Since $\mathcal P$ is an inductive partition we can find a sequence of convex bodies $H_1,\dots,H_l=K$  with the properties stated in Definition \ref{indpart}.
We argue by induction on $l$.
If $l=1$ the claim holds trivially. Assume now that the claim is true up to $l-1$.
If $H_l\in\mathcal P$ we can conclude as in the case $l=1$. Therefore we may assume $\exists j,k<l$ such that $H_j\cup H_k=K$ and $\interno(H_j\cap H_k)=\emptyset$.
As $H_j$ and $H_k$ are convex bodies, $u+I_{H_j}$ and $u+I_{H_k}$ belong to $\C^n$.
Moreover
$$(u+I_{H_j})\land (u+I_{H_k})=u+I_{H_j\cup H_k}=u+I_K,$$
while
$$(u+I_{H_j})\lor (u+I_{H_k})=u+I_{H_j\cap H_k}.$$
In particular, as $\interno(H_j\cap H_k)$ is empty, $\dim(\dom(u+I_{H_j\cap H_k}))<n$.
Hence, as $\mu$ is a simple valuation, we get
$$\mu(u+I_K)=\mu(u+I_{H_j})+\mu(u+I_{H_k}).$$
Now, if we set $\mathcal{P'}=\{P\in\mathcal P \,:\, P\sub H_j\}$ and $\mathcal{P''}=\{P\in \mathcal P \,:\, P\sub H_k\}$ and apply the inductive hypothesis to the just defined 
partitions we get
$$
\mu(u+I_{H_j})+\mu(u+I_{H_k})=\sum_{P\in\mathcal{P'}} \mu(u+I_P)+\sum_{P\in\mathcal{P''}}\mu(u+I_P)=\sum_{P\in\mathcal{P}}\mu(u+I_P).
$$ 
\end{proof}

\section{Characterization results I: simple valuations}\label{Characterization results I: simple valuations}

Our first characterization result is a converse of Proposition \ref{Valutazioni integrali semplici}. 

\begin{theorem}\label{charact-simple} Let $\mu\,:\,\C^n\to\R$ be a valuation with the following properties:
\begin{itemize}
\item[{\it i})] $\mu$ is rigid motion invariant;
\item[{\it ii})] $\mu$ is monotone decreasing;
\item[{\it iii})] $\mu$ is simple;
\item[{\it iv})] $\mu$ is m-continuous.
\end{itemize}
Then there exists a function $f\,:\,\R\to\R$, non-negative, decreasing, right-continuous and verifying the integrability condition
$$
\int_0^\infty t^nf(t)dt<\infty,
$$
such that for every $u\in\C^n$
\begin{equation}\label{charact-simple1}
\mu(u)=\int_{\dom(u)}f(u(x))dx.
\end{equation}
Equivalently
\begin{equation}\label{charact-simple2}
\mu(u)=\int_{\R}V_n(\cl(\{u< t\}))d\nu(t),
\end{equation}
where $\nu$ is the Radon measure related to $f$ by:
$$
f(t)=\nu((t,\infty))\quad\forall\,t\in\R.
$$ 
The function $f$ coincides with the geometric density
$f_n$ of $\mu$, determined by Proposition \ref{geometric-densities}.
\end{theorem}

Let us begin with some considerations preliminary to the proof. As $\mu$ is rigid motion invariant and monotone, its geometric densities
$f_i$, $i=0,\dots,n$ are defined (see Proposition \ref{geometric-densities}). On the other hand, by Proposition \ref{geometric densities simple valuations}, the only 
non-zero geometric density of $\mu$ is $f_n$. Recall that $f_n$ is a non-negative decreasing function defined on $\R$, moreover, as $\mu$ is m-continuous
 $f_n$ is right-continuous. Let $f\,:\,(-\infty,+\infty]\to[0,+\infty)$ be the extension of $f_n$, with the additional condition $f(\infty):=0$.
 
We will need the following Lemma.

\begin{lemma}\label{Riemann sums} Assume that $\mu$ is as in Theorem \ref{charact-simple} and let $f$ be the extension of its geometric density defined as above.
Let $K$ be a convex body and  
$$
{\mathcal P}=\{K_1,\dots,K_N\}
$$
be an inductive partition of $K$. 

Let $u\in\C^n$ such that $L=\dom(u)$ is a convex body, the restriction of $u$ to $L$ is continuous and 
$\interno (L)\sub K$. Then 
$$
\sum_{i=1}^N 
m_i\, V_n(K_i)\le
\mu(u+I_K)\le
\sum_{i=1}^N 
M_i\, V_n(K_i),
$$
where
$$
m_i=
\inf\{f(u(x))\,|\, x\in K_i\},\quad
M_i=
\sup\{f(u(x))\,|\, x\in K_i\}.
$$
\end{lemma}

\begin{proof} 
First we prove that $u$ attains a maximum and a minimum when restricted to $K_i$.

Since $u$ is lower semi-continuous and $K_i$ is compact and non-empty, we have $\inf_{K_i}u = \min_{K_i}u$.
Suppose now that there exists a point $x\in\K_i$ such that $u(x)=\infty$: in this case $\sup_{K_i} u=\infty$ is attained at $x$;
on the other hand, if  $K_i\sub L$, as the restriction of $u$ to $L$ (and thus to $K_i$) is continuous, $\sup_{K_i} u =\max_{K_i}u$.

Therefore, as $f$ is decreasing,
$$
M_i=f\left(\min_{K_i} u\right)\quad\mbox{and}\quad
m_i=f\left(\max_{K_i} u\right).
$$
Using the monotonicity of $\mu$ and the definition of geometric densities we obtain
\begin{eqnarray*}
\mu(u+I_{K_i})&\le&
\mu\left(\min_{K_i} u+I_{K_i}\right)=
M_i\,V_n(K_i),\\
\mu(u+I_{K_i})&\ge&
\mu\left(\max_{K_i} u+I_{K_i}\right)=
m_i\, V_n(K_i).
\end{eqnarray*}
Then 
\begin{equation}\label{proof Riemann sums 1}
m_i V_n(K_i)\le\mu(u+I_{K_i})\le
M_i V_n(K_i).
\end{equation}



By Lemma \ref{decomposition of simple valuations} we have that
$$
\mu(u+I_K)=
\sum_{i=1}^N \mu(u+I_{K_i}).
$$
This equality and \eqref{proof Riemann sums 1} conclude the proof.
\end{proof}

We will also need a well known theorem relating Riemann integrability and Lebesgue measure (known as Lebesgue-Vitali theorem). The proof of this result in 
the one-dimensional case can be found in standard texts of real analysis. 
The reader interested in the proof for general dimension may consult \cite{Allard}.

\begin{theorem}\label{VL} Let $g\,:\,\R^n\to\R$ be a bounded function which vanishes outside a compact set. Then $g$ is Riemann integrable 
in $\R^n$ if and only if the set of discontinuities of $g$ has Lebesgue measure zero.
\end{theorem}

\bigskip

\noindent{\em Proof of Theorem \ref{charact-simple}.}
Let us consider an arbitrary $u\in\C^n$ which, from now on will be fixed. If $\dim(\dom(u))<n$ then $f(u(x))=0$ for 
${\mathcal H}^n$-a.e. $x\in\R^n$; hence, as $\mu$ is simple,
$$
0=\mu(u)=\int_{\dom(u)} f(u(x))dx,
$$
i.e. the theorem is proven. Therefore in the remaining part of the proof we will assume that $\dim(\dom(u))=n$. 

Initially, we will assume that $\dom(u)=L$ is a convex body (with non-empty interior) and that the restriction of $u$ to
$L$ is continuous. This implies in particular that $g:=f\circ u$ has compact support.
We claim that the function $g$ is Riemann integrable on $\R^n$. This will follow
from Theorem \ref{VL} if we show that $g=f\circ u$ is bounded on $\R^n$ and the set of its discontinuities is a 
Lebesgue-null set.

It is easy to prove that $g$ is bounded, since it is non-negative by construction and, as $f$ is monotone decreasing,  
$\max_{\R^n}(g)=f(\min_{\R^n}(u))<\infty$. 

Since $f$ is monotone decreasing, the set of its discontinuities is countable: let us call it $\{\xi_i\}_{i\in\N}$ and set
$\xi_0=\infty$. 
We claim that the set of discontinuities of $g$ is contained in 
$$
\bigcup_{i\ge 0} \partial\{x\in\R^n \,:\, u(x)=\xi_i\}.
$$ 
Let $C\subseteq \R^n$ denote the set of points where $g$ is continuous. We therefore aim to prove the following:
$$
C^\complement \subseteq \bigcup_{i\ge 0} \partial(u^{-1}(\xi_i))=
\bigcup_{i\ge 0} \left(\cl(u^{-1}(\xi_i))\setminus \interno(u^{-1}(\xi_i))\right)=\bigcup_{i\ge 0} \left(\cl(u^{-1}(\xi_i))\cap (\interno(u^{-1}(\xi_i)))^\complement\right),
$$ 
which in turn is equivalent to 
$$
A:=\bigcap_{i\ge0}\left( (\cl(u^{-1}(\xi_i)))^\complement \cup \interno (u^{-1}(\xi_i))\right)\subseteq C.
$$
Let us take a fixed $x\in A$; for every choice of $i\ge 0$ there are two possibilities:
\begin{description}
  \item[(a)] $x\notin \cl(u^{-1}(\xi_i))$,
  \item[(b)] $x\in \interno(u^{-1}(\xi_i))$.
\end{description}

Suppose {\bf (b)} holds for two integers $i, j$; then 
$$
\left .\begin{array}{l l}
x\in \interno(u^{-1}(\xi_i))\subseteq u^{-1}(\xi_i) \implies u(x)=\xi_i \\
x\in \interno(u^{-1}(\xi_j))\subseteq u^{-1}(\xi_j) \implies u(x)=\xi_j 
\end{array}\right\} \text{ then }i=j.
$$
Which means that {\bf (b)} can happen at most once for every choice of $x$. Let us prove $x\in C$ in case {\bf (b)} never holds.
As $x\notin \cl(u^{-1}(\xi_0))$,  and $\xi_0=\infty$, $x$ is an interior point of the domain of $u$, so that $u$ is continuous at $x$.
Moreover, for every $i\ge 0$ we have $x\notin \cl(u^{-1}(\xi_i))$ which implies $x\notin u^{-1}(\xi_i)$, i.e. $f$ is continuous at $u(x)$. It follows that
$g$ is continuous at $x$. It remains to prove $x\in C$ in case {\bf (b)} holds for a specific $j\ge 0$.
Since $x\in \interno(u^{-1}(\xi_j))$ there exists a neighborhood $B$ of $x$ such that $B\subseteq \interno(u^{-1}(\xi_j))\subseteq u^{-1}(\xi_j)$, 
which means $u(B)=\{\xi_j\}$. Thus $u$, and also $g$, is constant on a neighborhood of $x$ and hence continuous at $x$.

Having finally proven $A\subseteq C$ it remains to show that $A^\complement$ is a Lebesgue-null set.
Since $\partial u^{-1}(\xi_i)$ is the boundary of a convex body (possibly being empty) for every $i$, all these sets must be null sets, and so is their countable union (and therefore $C$ 
itself is a null set, being a subset of a null set). 

Now let $R$ be a closed rectangle such that $\interno(L)\subset R$; in particular $g$ vanishes in the complement of $R$. As $f$ is Riemann integrable, for every
$\epsilon>0$ there exists a rectangular partition (see section \ref{Rectangular partitions}) of $R$
$$
{\mathcal P}=\{R_1,\dots, R_N\}
$$
such that 
$$
\sum_{i=1}^N\left(
\sup_{R_i} g-\inf_{R_i} g
\right)V_n(R_i)\le\epsilon,
$$
and, clearly,
$$
\sum_{i=1}^N\inf_{R_i} g\, V_n(R_i)\le
\int_{\R^n} g\,dx\le
\sum_{i=1}^N\sup_{R_i} g\, V_n(R_i).
$$
Without loss of generality we could assume that $\mathcal P$ is an inductive partition and thus, by Lemma \ref{Riemann sums} we have
$$
\sum_{i=1}^N\inf_{R_i} g\, V_n(R_i)\le
\mu(u)\le
\sum_{i=1}^N\sup_{R_i} g\, V_n(R_i)
$$
(note that $u+I_R=u$, as $R$ contains the domain of $u$). Hence
$$
\left|
\mu(u)-\int_{\R^n}f(u(x))\,dx
\right|\le \epsilon
$$ 
and, as $\epsilon$ is arbitrary,
$$
\mu(u)=\int_{\R^n}f(u(x))\,dx,
$$ 
i.e. \eqref{charact-simple1} is true for every $u\in\C^n$ such that the domain of $u$ is a compact convex set and $u$ is continuous on its domain.
To prove the same equality for a general $u$, let $L_i$, $i\in\N$, be an increasing sequence of convex bodies such that
$$
\bigcup_{i=1}^\infty L_i=\dom(u),
$$
and consider the sequence of functions $u_i$, $i\in\N$, defined by $u_i=u+I_{L_i}$. This is a decreasing sequence of elements of $\C^n$ converging point-wise to
$u$ in the interior of $\dom(u)$. By the m-continuity of $\mu$ we have
$$
\mu(u)=\lim_{i\to\infty}\mu(u_i)=
\lim\int_{\R^n} f(u_i)\,dx
$$
where, in the second equality, we have used the first part of the proof. On the other hand the sequence of functions $f\circ u_i$, $i\in\N$, is increasing and 
converges point-wise to $f\circ u$ in $\R^n$. Hence, by the monotone convergence theorem,
$$
\lim_{i\to\infty}\int_{\R^n} f(u_i)\,dx=
\int_{\R^n} f(u)\,dx.
$$
The proof of \eqref{charact-simple1} is complete. As for \eqref{charact-simple2}, it follows from Proposition \ref{i=n}.
\begin{flushright}
$\qed$
\end{flushright}

\begin{remark}\label{osservazione1} 
{\em
It is clear from the previous proof that the representation formula \eqref{charact-simple1} of Theorem \ref{charact-simple} remains valid 
for those functions $u\in\C^n$, such that: $\dom(u)=L\in\K^n$ and the restriction of $u$ to $L$ is continuous,
even if we drop the assumption of m-continuity of $\mu$.}
\end{remark}

\section{Characterization results II: homogeneous valuations}\label{Characterization results II: homogeneous valuations}. 

\subsection{Part one: $n$-homogeneous valuations} 

The following result is a direct consequence of Theorem \ref{charact-simple} and Proposition \ref{k-homogeneous implies k-simple}.

\begin{theorem}\label{charact-n hom} Let $\mu\,:\,\C^n\to\R$ be a valuation with the following properties:
\begin{itemize}
\item $\mu$ is rigid motion invariant;
\item $\mu$ is monotone decreasing;
\item $\mu$ is $n$-homogeneous;
\item $\mu$ is m-continuous.
\end{itemize}
Then there exists a function $f\,:\,\R\to\R$, coinciding with the geometric density
$f_n$ of $\mu$, non-negative, decreasing, right-continuous and verifying the integrability condition
$$
\int_0^\infty t^nf(t)dt,
$$
such that for every $u\in\C^n$
$$
\mu(u)=\int_{\dom(u)}f(u(x))dx,
$$
or equivalently
$$
\mu(u)=\int_{\R}V_n(\cl(\{u<t\}))d\nu(t),
$$
where $\nu$ is the Radon measure related to $f$ by
$$
f(t)=\int_{(t,\infty)}d\nu(s)\quad\forall\, t\in\R.
$$ 
\end{theorem}

\begin{definition}{\bf(Extensions and restrictions of convex functions).}
Let $k<n$. Let $u\in\C^k$. We can now extend $u$ to the whole $\rn$ in a canonical way by assigning the value $\infty$ where $u$ was otherwise undefined:
$$
\ext{u}{n}(x) 
=\ext{u}{n}(x_1,\dots,x_k,x_{k+1},\dots,x_n)=
 \begin{cases} 
      \hfill u(x_1,\dots,x_k)    \hfill & \text{ if $x_{k+1}=\dots=x_n=0$} \\
      \hfill \infty \hfill & \text{otherwise} \\
  \end{cases}.$$  
If $u\in\C^k$, then, it can be shown that $\ext{u}{n}\in\cn$.
On the other hand, the so-called restriction of a convex function $u\in\C^n$ can be defined in the following way:
$$ \restr{u}{k} (x_1,\dots, x_k)=u(x_1,\dots,x_k,0,\dots,0).$$
It is immediate to show that $\restr{u}{k}$ belongs to $\C^k$ for every choice of $u\in\C^n$.
\end{definition}
\begin{definition}{\bf (Restrictions of valuations).}\label{restrval}
Let $k<n$ as above. Let $\mu$ be a real valuation on $\C^n$, then we can define the restriction of $\mu$ to $\C^k$ as
$$\restr{\mu}{k}(u)=\mu(\ext{u}{n}) \;\;\forall u\in\C^k.$$
\end{definition}
It is easy to verify that $\restr{\mu}{k}$ defined as above is a valuation on $\C^k$. Moreover, the valuation $\restr{\mu}{k}$ inherits the following properties from $\mu$: rigid motion invariance, monotonicity, m-continuity and homogeneity.
Let us now consider a valuation $\mu$ on $\C^n$ and a convex function $u\in\C^n$ such that 
\begin{equation}\label{dompicc}
\dom(u)\sub \{ (x_1,\dots, x_n)\in \R^n \,:\, x_{k+1}=\dots=x_n=0\}.
\end{equation}
Under these assumptions we have that 
\begin{equation}\label{restrest}
\mu(u)=\restr{\mu}{k}(\restr{u}{k}).
\end{equation}
The previous equality is an immediate consequence of Definition \ref{restrval} and of the following consideration:
$$\ext{(\restr{u}{k})}{n}=u$$ 
for every $u\in\C^n$ which satisfies \eqref{dompicc}.
Restricted valuations also share geometric densities up to the suitable dimension. To be more precise, if $f_0,\dots,f_n$ are the geometric densities of $\mu$, then $f_0,\dots,f_k$ are the geomtric densities of $\restr{\mu}{k}$.
To prove this, let $t$ be a real number and $H$ an arbitrary convex body in $\K^k$. 
Then, 
$$K\times P \in \K^n$$
where $P=\{(0,\dots,0)\}\in\K^{n-k}$.
We have 
$$\dom(t+I_{K\times P})=K\times P,$$
in other words, the function $t+I_{K\times P}$ satisfies \eqref{dompicc}.
We deduce that \eqref{restrest} holds for $u=t+I_{K\times P}$, thus
$$
\restr{\mu}{k}\left( \restr{(t+I_{K\times P})}{k}\right) = \mu(t+I_{K\times P})= \sum_{i=0}^n f_i(t)V_i(K\times P)=\sum_{i=0}^k f_i(t)V_i(K).
$$ 
A simple calculation yields
$$ \restr{(t+I_{K\times P})}{k} = t+ I_K\in\C^k.$$
Therefore 
$$
\restr{\mu}{k}(t+I_K)=\sum_{i=0}^k f_i(t)V_i(K),
$$
we conclude by the arbitrariness of $t$ and $K$.

The following corollary of Theorem \ref{charact-n hom} will be important in the sequel.

\begin{corollary}\label{corollary k-hom} Let $\mu\,:\,\C^n\to\R$ be a valuation with the following properties:
\begin{itemize}
\item $\mu$ is rigid motion invariant;
\item $\mu$ is monotone decreasing;
\item $\mu$ is $k$-homogeneous, for some $k\in\{0,1,\dots,n\}$;
\item $\mu$ is m-continuous.
\end{itemize}
Let $f_k$ denote the $k$-th geometric density of $\mu$. Then $f_k$ verifies the integrability condition
$$
\int_0^\infty f_k(t)^{k-1}dt<\infty.
$$
Moreover, for every $u\in\C^n$ such that $\dim(\dom(u))\le k$
$$
\mu(u)=\int_{\R} V_k(\cl(\{u < t\}))d\nu(t)
$$
where $\nu$ is a Radon measure on $\R$ and $f_k$ and $\nu$ are related by the identity
$$
f_k(t)=\int_{(t,\infty)} d\nu(s)\quad\forall\, t\in\R.
$$
\end{corollary}

\begin{proof} Starting from $\mu$ we define its restriction to $\C^k$, $\restr{\mu}{k}$. 

As remarked, $\restr{\mu}{k}$ is a valuation on $\C^k$ with the following properties: it is rigid motion invariant, monotone decreasing, m-continuous
and $k$-homogeneous.
 Denote by $g_i$, $i\in\{0,1,\dots,k\}$ its geometric densities; then $g_0\equiv\dots g_{k-1}\equiv0$ and
$g_k=f_k$. By Theorem \ref{charact-n hom} we have that $f_k$ verifies the claimed integrability condition and
\begin{equation}\label{proof4}
\restr{\mu}{k}(v)=\int_{\R}V_k(\cl(\{v<t\}))d\nu(t)
\end{equation}
for every $v\in\C^k$, where $\nu$ and $f_k$ are related as usual by $f_k(t)=\nu((t,\infty))$. Now let $u\in\C^n$ be such that 
$\dim(\dom(u))\le k$; we want to compute
$\mu(u)$. As $\mu$ is rigid motion invariant, without loss of generality we may assume that 
$$
\dom(u)\sub\{(x_1,\dots,x_n)\in\R^n\,:\, x_{k+1}=\dots=x_n=0\}.
$$
Then $\mu(u)=\restr{\mu}{k}(\restr{u}{k})$. 
If we set $P=\{(0,\dots,0)\}\in\K^{n-k}$ as before, then it is simple to verify that
$$\{u<t\}= \{\restr{u}{k}<t\}\times P \implies \cl(\{u<t\})=\cl(\{\restr{u}{k}<t\})\times P,$$
so that $V_k(\cl(\{u<t\}))=V_k( \cl(\{\restr{u}{k}<t\}))$.
The claimed representation formula for $\mu(u)$ follows from the previous considerations and \eqref{proof4} specialized to the case $v=\restr{u}{k}$.
\end{proof}

\subsection{Part two: the general case}

\begin{theorem}\label{characterization of k-homogeneous valuations} Let $\mu\,:\,\C^n\to\R$ be a valuation with the following properties:
\begin{itemize}
\item $\mu$ is rigid motion invariant;
\item $\mu$ is monotone decreasing;
\item $\mu$ is m-continuous;
\item $\mu$ is $k$-homogeneous, for some $k\in\{0,\dots,n\}$.
\end{itemize}
Then there exists a Radon measure $\nu$ defined on $\R$, verifying the integrability assumption
$$
\int_{\R}t^{k}d\nu(t)<\infty,
$$
such that
$$
\mu(u)=\int_{\R}V_{k}(\cl(\{u< t\}))d\nu(t),
$$
for every $u\in\C^n$. Moreover, the measure $\nu$ is determined by the unique non-vanishing geometric density $f_k$ of $\mu$ as follows:
$$
f_k(t)=\int_{(t,\infty)} d\nu(s)\quad\forall\, t\in\R.
$$
\end{theorem}

The rest of this section is devoted to the proof of this result; throughout, $\mu$ will be a valuation with the properties indicated in the previous theorem.
Note that the validity of the Theorem for $k=n$ is established by Theorem \ref{charact-n hom}. 

By Proposition \ref{geometric-densities} we may assign to $\mu$ its geometric densities $f_j$, $j=0,\dots,n$. 
By homogeneity we have 
$f_j\equiv 0$ for every $j\ne k$. In other words, the only density which can be non-identically 
zero is $f_{k}$. For simplicity we will call this function $f$. By the properties of $\mu$, this is a non-negative decreasing function; moreover,
as $\mu$ is m-continuous, $f$ is right-continuous on $\R$ and, by Corollary \ref{corollary k-hom}, it verifies the integrability condition
$$
\int_0^\infty t^{k-1}f(t)dt<\infty.
$$

We proceed by induction on the dimension $n$. Let us then start from the case $n=1$. As the theorem is already proven for $k=n=1$ we only need to consider the case
$k=0$; but this follows from Proposition \ref{characterization of zero homogeneous valuations}: Theorem \ref{characterization of k-homogeneous valuations} is proven in
dimension $n=1$. 

\medskip

To continue with the induction argument, we assume that the theorem holds up to dimension $(n-1)$ and we are going to prove it in the $n$-dimensional case.  We may assume
that $1\le k\le n-1$. 

In the next part of the proof we will assume, in addition to the above properties, that the only non-zero density $f$ of $\mu$ is smooth: $f\in C^\infty(\R)$. 
As in the previous sections, we introduce the Radon measure $\nu$ related to $f$ by the identity
$$
f(t)=\int_{(t,\infty)} d\nu(s)\quad\forall\, t\in\R.
$$

Based on $\nu$, we construct an auxiliary valuation $\mu_a\,:\,\C^n\to\R$ defined as follows
$$
\mu_a(u)=\int_{\R} V_{k}(\cl(\{u<t\}))d\nu(t).
$$
By the results of section \ref{sec geometric densities}, $\mu$ is a well defined valuation and it is rigid motion invariant, decreasing, $k$-homogeneous
and m-continuous. Moreover, its geometric density of order $k$ is precisely $f$, i.e. the same as $\mu$.

We also set
$$
\mu_r=\mu_a-\mu.
$$
The idea is to prove that $\mu_r$ is identically zero. Note that $\mu_r$ inherits most of the properties of $\mu$ and $\mu_a$: it is a valuation,
rigid motion invariant, $k$-homogeneous and m-continuous. We cannot infer in general that $\mu_r$ is monotone. 

\bigskip

\noindent
{\bf Claim 1.} {\em The valuation $\mu_r$ ``vanishes horizontally'', i.e. for every convex body $K\in\K^n$ and every $t\in\R$ we have
$$
\mu_r(t+I_K)=0.
$$}

The proof is a straightforward consequence of the fact that $\mu$ and $\mu_a$ have the same geometric densities.

\bigskip

\noindent{\bf Claim 2.} {\em The valuation $\mu_r$ is simple.}

\begin{proof}
Let $u\in\C^n$ be a convex function whose domain has dimension strictly less than $n$. As $\mu_r$ is rigid motion invariant, we might assume without loss of generality that 
$$u\sub \{(x_1,\dots, x_n)\in\R^n \,:\, x_n=0\}.$$ 
As remarked after Definition \ref{restrval}, $\mu(u)=\restr{\mu}{n-1}(\restr{u}{n-1})$.
By the induction hypothesis we have
\begin{equation}\label{questa}
\restr{\mu}{n-1}(\restr{u}{n-1})=\int_{\R} V_k(\cl(\{\restr{u}{n-1}<t\}))d\nu(t).
\end{equation}
As $\{u<t\}=\{\restr{u}{n-1}<t\}\times \{0\}$, \eqref{questa} can be rewritten as
$$ \mu(u)=\restr{\mu}{n-1}(\restr{u}{n-1})=\int_{\R} V_k(\cl(\{u<t\}))d\nu(t)=\mu_a(u). $$
Therefore $\mu_r(u)=0$, we conclude that $\mu_r$ is simple as claimed.  
\end{proof}

\bigskip

We will now introduce a construction which is going to help us evaluate a valuation on piece-wise linear functions.
We fix 
$$
\mbox{$e\in\R^n$ s.t. $|e|=1$,}\quad p\ge0,\quad V=pe.
$$
Let $\mu_0$ be a valuation on $\C^n$ and consider the linear function $w\,:\,\R^n\to\R$ defined by
$$
w(x)=(x,V),\quad x\in\R^n.
$$
Then we define a mapping on the family of convex bodies of $\R^n$, $\sigma_{\mu_0,V}\,:\,\K^n\to\R$, as follows
$$
\sigma_{\mu_0,V}(K)=\mu_0(w+I_K)\quad\forall\, K\in\K^n.
$$
It is easy to check that $\sigma_{\mu_0,V}$ is a valuation.

From now on throughout this paper, we will consider the previous construction specialized to valuations which are rigid motion invariant, hence we will assume
without loss of generality that $V=pe_n=(0,\dots,0,p)\in\R^n$, $p\ge 0$, and consequently that $w(x)=w_p(x)=p(x,e_n)=px_n$ for all $x\in\R^n$. 
Moreover, for the sake of brevity, we will introduce the following simplified notation:
$$
\sigma:=\sigma_{\mu,V},\quad
\sigma_a:=\sigma_{\mu_a,V},\quad
\sigma_r:=\sigma_{\mu_r,V}=\sigma-\sigma_a.
$$

The following claim collects some of the properties of $\sigma_r$ that will be
used in the sequel.

\bigskip

\noindent{\bf Claim 3.} {\em $\sigma_r$ has the following properties:
\begin{itemize}
\item[1.] it is a valuation on $\K^n$;
\item[2.] it is simple;
\item[3.] it is invariant with respect to every rigid motion $T$ of $\R^n$ 
such that
\begin{equation}\label{claim3 I}
T(x)=T(x_1,\dots,x_n)=
(T'(x_1,\dots,x_{n-1}),x_n)\quad\forall\,x\in\R^n,
\end{equation}
where $T'$ is a rigid motion of $\R^{n-1}$.
\end{itemize}}

\begin{proof} Let $K,L\in\K^n$ be such that $K\cup L\in\K^n$. Then 
$$
I_{K\cup L}=I_K\wedge I_L,\quad
I_{K\cap L}=I_K\vee I_L.
$$
These relations remain valid if we add $w$ as follows
$$
w+I_{K\cup L}=(w+I_K)\wedge(w+I_L),\quad
w+I_{K\cap L}=(w+I_K)\vee (w+I_L).
$$
Using the valuation property of $\mu_r$ we easily deduce the valuation property for $\sigma_r$. 
Moreover 
$$
\sigma_{r}(\emptyset)=\mu_r(w+I_\emptyset)=\mu_r(\boldsymbol{\infty})=0.
$$
We conclude that $\sigma_r$ is a valuation.

If $K\in\K^n$ has no interior point, the domain of $I_K$, and consequently that of $w+I_K$, have the same property.
Then,
$$
\sigma_r(K)=\mu_r(w+I_K)=0.
$$

Next we prove 3. Let $T$ be a rigid motion of $\R^n$ of the form \eqref{claim3 I}, and let $K\in\K^n$. Then
\begin{eqnarray*}
(w+I_{T(K)})(x)&=&px_n+I_{T(K)}(x)=px_n+I_K(T'^{-1}(x_1,\dots,x_{n-1}),x_n)\\
&=&(w+I_K)(T'^{-1}(x_1,\dots,x_{n-1}),x_n)=(w+I_K)(T^{-1}(x)).
\end{eqnarray*}
Therefore
$$
\sigma_r(T(K))=\mu_r(w+I_{T(K)})=\mu_r((w+I_K)\circ T^{-1})
=\mu_r(w+I_K)=\sigma_r(K),
$$
where we have used the invariance if $\mu_r$.
\end{proof}

We anticipate that the following step is one of the most delicate in the proof.

\bigskip

\noindent{\bf Claim 4.} 
{\em The valuation $\sigma_r$ is non-negative, i.e.}
$$
\sigma_{r}(K)\ge0\quad\forall\, K\in\K^n.
$$

\begin{proof} 

We first treat the easier case $p=0$, which leads to $w\equiv 0$ so that
$$
\sigma (K)=\mu(I_K)=f(0)V_{k}(K)=\sigma_a(K)\,\Rightarrow
\,\sigma_r(K)=\sigma_a(K)-\sigma(K)=0,
$$
where we have used Claim 1. Next we assume $p>0$. 
Given two real numbers $\alpha,\beta$ with $\alpha\le\beta$, we define the strip:
$$
S[\alpha,\beta]:=\{(x_1,\dots,x_n)\in\R^n\,:\,\alpha\le x_n\le\beta\}.
$$
Let $K\in\K^n$ and let $y_m$ and $y_M$ be such that the hyperplanes of equations
$x_n=y_m$ and $x_n=y_M$ are the supporting hyperplanes to $K$ with outer unit normals
$-e_n$ and $e_n$ respectively. In other words, 
$$
K\subset S[y_m,y_M]
$$
and $S[y_m,y_M]$ is the intersection of all possible sets of the form $S[\alpha, \beta]$ containing $K$. We define a function 
$\phi\,:\,[y_m,y_M]\to\R$: 
$$
\phi(y)=\sigma_r(K\cap S[y_m,y]).
$$
Let $y\in(y_m,y_M)$ and $h\ge0$ be sufficiently small so that $y+h\le y_M$.
As, trivially,
$$
S[y_m,y]\cup S[y,y+h]=S[y_m,y+h],\quad
S[y_m,y]\cap S[y,y+h]=S[y,y],
$$
using the valuation property of $\sigma_r$ and the fact that it is simple, we get
\begin{eqnarray*}
\phi(y+h)-\phi(y)&=&\sigma_r(K\cap S[y,y+h])\\
&=&\sigma_a(K\cap S[y,y+h])-\sigma (K\cap S[y,y+h]).
\end{eqnarray*}
Next we use the monotonicity of $\mu$. Note that
$$
\min_{K\cap S[y,y+h]}w=py,\quad
\max_{K\cap S[y,y+h]}w=p(y+h).
$$
Hence
\begin{eqnarray*}
\sigma(K\cap S[y,y+h])&=&\mu(w+I_{K\cap S[y,y+h]})
\le\mu(py+I_{K\cap S[y,y+h]})\\
&=&f(py)\, V_{k}(K\cap S[y,y+h]).
\end{eqnarray*}
And similarly 
$$
\sigma(K\cap S[y,y+h])\ge f(p(y+h))\,V_{k}(K\cap S[y,y+h]).
$$
On the other hand
\begin{eqnarray*}
\sigma_a(K\cap S[y,y+h])&=&\mu_a(w+I_{K\cap S[y,y+h]})\\
&=&-\int_{\R}V_{k}(\cl(\{w+I_{K\cap S[y,y+h]}<t\}))f'(t)dt,
\end{eqnarray*}
where we have used the assumption that $f$ is smooth. Now
$$
\cl(\{w+I_{K\cap S[y,y+h]}< t\})=
\left\{
\begin{array}{lllll}
\emptyset\quad&\mbox{if $t\le py$,}\\
\\
K\cap S[y,t/p]\quad&\mbox{if $py< t< p(y+h)$,}\\
\\
K\cap S[y,y+h]\quad&\mbox{if $t\ge p(y+h)$.}
\end{array}
\right.
$$

Hence
\begin{eqnarray*}
\sigma_a(K\cap S[y,y+h])&=&
-\int_{py}^{p(y+h)}V_{k}(K\cap S[y,t/p])f'(t)dt
-V_{k}(K\cap S[y,y+h])\int_{p(y+h)}^\infty f'(t)dt\\
&=&-p\int_{y}^{y+h}V_{k}(K\cap S[y,s]) f'(ps)ds+
f(p(y+h))V_{k}(K\cap S[y,y+h]),
\end{eqnarray*}
where, for the second term we have used the equality, due to the condition $k\ge 1$ and the integrability condition on $f$,
$$
\lim_{t\to\infty}f(t)=0.
$$
Consequently we have the following bounds:
\begin{equation}\label{proof5}
\phi(y+h)-\phi(y)\le -p\int_y^{y+h}V_{k}(K\cap S[y,y+h])f'(ps)ds,
\end{equation}
and
\begin{eqnarray}\label{proof6}
\phi(y+h)-\phi(y)&\ge&V_{k}(K\cap S[y,y+h])\left(
f(p(y+h))-f(py)\right)\\
&-&p\int_y^{y+h}V_{k}(K\cap S[y,y+h])f'(ps)ds=0.\nonumber
\end{eqnarray}

Note that the function 
$$
\tau\,\mapsto\,V_{k}(K\cap S[y,y+\tau])
$$
is Lipschitz continuous in a neighborhood of $\tau=0$ (indeed, as already remarked before, its $1/k$ power is concave in
$[y_m,y_M]$) and, by monotonicity of intrinsic volumes, it is bounded by $V_{k}(K)$, i.e. a constant independent of $y$ and $h$. 
Then, as $f$ is smooth, it follows from \eqref{proof5} and \eqref{proof6} that $\phi$ is Lipschitz continuous in $[y_m,y_M]$; in particular
\eqref{proof6} implies that
$$
\phi'(y)\ge 0
$$
for every $y$ for which $\phi'$ is defined. As 
$$
\phi(y_m)=\sigma_r(K\cap S[y_m,y_m])=0
$$ 
(recall that $\sigma_r$ is simple), we have that
$$
\sigma_r(K)=\phi(y_M)=\phi(y_m)+\int_{y_m}^{y_M}\phi'(t)dt
\ge0.
$$ 
The proof is complete. 
\end{proof}

\bigskip

For the sequel we will need the following result (which could be well-known in the theory of valuations on convex bodies).

\begin{lemma}\label{positive implies monotone} Let $\sigma\,:\,\K^n\to\R$ be a valuation which is
non-negative and simple. Then $\sigma$ is monotone increasing on the class of polytopes, i.e. for every $P$ and $Q$
polytopes in $\R^n$ such that $P\sub Q$ we have
$$
\sigma(P)\le \sigma(Q).
$$
\end{lemma}
\begin{proof} Let $P$ and $Q$ be polytopes such that $P\sub Q$; let $\mathcal F$ be a family of hyperplanes in $\R^n$, defined as follows:
$$
{\mathcal F}(P,Q)=\left\{\mbox{$H$ is a hyperplane containing a facet of $P$ and $H\cap\interno(Q)\ne\emptyset$}\right\},
$$
and let $N(P,Q)\ge0$ be the cardinality of ${\mathcal F}(P,Q)$. We will prove that
$$
\sigma(P)\le\sigma(Q),
$$
by induction on $N(P,Q)$. If $N(P,Q)=0$ we have that $P=Q$ so that there is nothing to prove.
Assume that the claim is true up to $(n-1)$, for some $n\in\N$, and that $N(P,Q)=n$. Let 
$H\in{\mathcal F}$ and let $H^+$ and $H^-$ be the closed half-spaces determined by $H$. We may assume that
$P\sub H^+$. Let $Q^+=Q\cap H^+$ and $Q^-=Q\cap H^-$ (which are still polytopes); as
$$
Q=Q^+\cup Q^-\quad\mbox{and}\quad
Q^+\cap Q^-\subset H,
$$
and as $\sigma$ is simple and non-negative
$$
\sigma(Q)=\sigma(Q^+)+\sigma(Q^-)\ge\sigma(Q^+).
$$
On the other hand $Q^+\supseteq P$ and 
$$
{\mathcal F}(P, Q^+)\subset {\mathcal F}(P,Q),
$$
in particular $N(P,Q^+)<N(P,Q)$ so that, by the induction assumption,
$$
\sigma(P)\le\sigma(Q^+)\le\sigma(Q).
$$
\end{proof}

\bigskip

\noindent
{\bf Claim 5.} 
{\em The valuation $\sigma_r$ is monotone increasing.}
\begin{proof} Let $K,L\in\K^n$ be such that $K\sub L$; there exist two sequences of polytopes $P_i$, $Q_i$, $i\in\N$, with the following properties:
\begin{enumerate}
\item they are increasing with respect to set inclusion;
\item $P_i\to K$ and $Q_i\to L$ as $i$ tends to infinity, in the Hausdorff metric;
\item $P_i\sub Q_i$ for every $i\in\N$.
\end{enumerate}
In particular $\sigma_r (P_i)\le \sigma_r (Q_i)$. Now, recalling the definition of $\sigma_{r}$ we have that
\begin{eqnarray*}
\sigma_{r}(P_i)&=&\sigma_{a}(P_i)-\sigma(P_i)
=\mu_a(u_i)-\mu(u_i)
\end{eqnarray*}
where
$$
u_i(x)=I_{P_i}(x)+(x,V),\quad x\in\R^n.
$$
Note that $u_i$ is a decreasing sequence, and it converges point-wise to $u\,:\,\R^n\,\to\R$ defined by
$$
u(x)=I_K(x)+(x,V)
$$
in the relative interior of $K$. As $\mu_a$ and $\mu$ are m-continuous we have
$$
\sigma_{r}(P_i)=
\mu_r(u_i)=\mu_a(u_i)-\mu(u_i)\to\mu_a(u)-\mu(u)=\mu_r(u)=\sigma_{r}(K).
$$
In a similar way we can prove that
$$
\lim_{i\to\infty}\sigma_{r}(Q_i)=\sigma_{r}(L).
$$
Since, as already pointed out, $\sigma_r(P_i)\le \sigma_r(Q_i)$ for all $i\in\N$, passing to the limit for $i\to\infty$ yields the claimed 
$\sigma_r(K)\le \sigma_r(L).$ 
\end{proof}

\bigskip

Let us make a further step to investigate the behavior of $\mu_r$ on restrictions of linear functions. 
Given a function $u\in\C^{n-1}$, we may consider the set
$$
\epi(u)=\{(x',y)\in \R^{n-1}\times\R\,:\,y\ge u(x')\}.
$$

\bigskip
\noindent{\bf Claim 6.} {\em Let $p>0$. For every $u\in\C^{n-1}$ the function $w+I_{\epi(u)}$ belongs to $\C^n$.}

\begin{proof}
The set $\epi(u)$ is convex and closed (by the semi-continuity of $u$). Hence the function
$v=w+I_{\epi(u)}$ is lower semi-continuous and convex. Let $x_i$, $i\in\N$, be a sequence in $\R^n$ such that 
$$
\lim_{i\to\infty}|x_i|=\infty.
$$
From any subsequence of $x_i$ we may extract a further subsequence (let us call it $\bar x_i$) such that either
$\bar x_i\in\epi(u)$ for every $i$ or $\bar x_i\in\R^n\setminus\epi(u)$ for every $i$.  In the second case we have $v(\bar x_i)=\infty$ for every $i$.
In the first case we have, setting $\bar x_i=(\bar x_i',\bar y_i)\in V^\perp\times\R$, there exists constants $a>0$, $b\in\R$ such that
$$
\bar y_i\ge u(\bar x_i)\ge a|\bar x_i'|+b\quad\forall\, i\in\N
$$
(see Proposition \ref{limitazione-uniforme}). As $|\bar x_i|$ is unbounded, we must have that $\bar y_i$ is not bounded from above and, up to extracting a further 
subsequence
we may assume that 
$$
\lim_{i\to\infty} \bar y_i=\infty.
$$
This implies that $v(\bar x_i)=p y_i$ tends to infinity as well. Hence from any subsequence $x_i$ such that $|x_i|\to\infty$ we may extract a subsequence 
$\hat x_i$ such that $v(\hat x_i)$ tends to infinity. Hence
$$
\lim_{|x|\to\infty} v(x)=\infty
$$
and we conclude that $v\in\C^n$.
\end{proof}

Let us define $\bar\mu\,:\C^{n-1}\to\R$ by
$$
\bar\mu(u)=\mu_r(w+I_{\epi(u)}).
$$

\bigskip

\noindent
{\bf Claim 7.} {\em Let $p>0$. The application $\bar\mu$ has the following properties:
\begin{itemize}
\item[1.] it is a rigid motion invariant valuation;
\item[2.] it is simple;
\item[3.] it is monotone decreasing.
\end{itemize}}

\begin{proof} 
We will denote a point in $\R^n$ by $(x',y)$, with $x'\in\R^{n-1}$ and $y\in\R$. For $u\in\C^{n-1}$ and $t\in\R$ set
$$
\epi_t(u)=\{(x',y)\in\R^{n-1}\times\R\,:\, u(x')\le y\le t\}=\epi(u)\cap\{(x',y)\,:\, y\le t\}.
$$
By the m-continuity of $\mu_r$,
$$
\bar\mu(u)=\lim_{t\to\infty}\mu_r(w+I_{\epi_t(u)})
=\lim_{t\to\infty}\sigma_{r}(\epi_t(u)).
$$
We will see that properties 1 - 3 follow easily from this characterization of $\bar\mu$ and Claim 3.
Assume that $T$ is a rigid motion of $\R^{n-1}$. Define
$$
\bar T\,:\,\R^n\to\R^n,\quad
\bar T(x',y)=(T(x'),y).
$$
$\bar T$ is a rigid motion of $\R^n$ and it verifies $(\bar T(x),V)=(x,V)$ for every $x\in\R^n$. Then, by item 3 in Claim 3,
$$
\sigma_{r}(\epi_t(u))=\sigma_{r}(\bar T(\epi_t(u))).
$$
On the other hand
$$
\bar T(\epi_t(u))=\epi_t(u\circ T^{-1}).
$$
Replacing this equality in the previous one, and letting $t\to\infty$, we get
$$
\bar\mu(u)=\bar\mu(u\circ T^{-1}),
$$
which proves that $\bar\mu$ is rigid motion invariant.

To prove that $\bar\mu$ is simple, let $u\in\C^{n-1}$ be such that $\dim(\dom(u))\le(n-2)$. Then
$\dim(\epi(u))\le(n-1)$ and
$$
\dim(\epi_t(u))\le n-1\quad\forall\,t\in\R.
$$
As $\mu_r$ is simple, $\sigma_{r}$ is simple, by Claim 3. Hence
$$
\sigma_{r}(\epi_t(u))=0\quad\forall\,t.
$$
Letting $t$ tend to infinity we get $\bar\mu(u)=0$. 

As for monotonicity, if $u$ and $v$ belong to $\C^{n-1}$ and are such that $u\le v$ in $\R^{n-1}$, then 
$$
\epi(u)\supseteq\epi(v)\;\Rightarrow\;
\epi_t(u)\supseteq\epi_t(v)\quad\forall\,t\in\R.
$$
As $\sigma_{r}$ is monotone increasing we get
$$
\sigma_{r}(\epi_t(u))\ge\sigma_{r}(\epi_t(v))\quad\forall\, t\in\R.
$$
The conclusion follows letting $t$ to $\infty$.
\end{proof}

\bigskip

\noindent
{\bf Claim 8.} {\em Let $p>0$. There exists a function $\bar f\,:\,\R\times(0,\infty)$, $\bar f=
\bar f(t,p)$, such that 
\begin{equation}\label{Claim 8}
\bar\mu(u)=\int_{\dom(u)}\bar f(u(x'),p)dx'
\end{equation}
for every $u\in\C^{n-1}$ such that $\dom(u)\in\K^n$ and the restriction of $u$ to $\dom(u)$ is continuous
(here $dx'$ denotes the usual integration in $\R^{n-1}$).} 
\begin{proof}
By Claim 7 we may apply Theorem \ref{charact-simple} and subsequent Remark \ref{osservazione1} to deduce \eqref{Claim 8}. 
\end{proof}

\medskip

Given $K\in\K^{n-1}$ and $t_1,t_2\in\R$, with $t_1\le t_2$,
we consider the cylinder:
$$
K\times [t_1,t_2]\in\K^n.
$$
Evaluating $\sigma_r$ on cylinders is a crucial step, as we will see in the following claim.
\bigskip

\noindent
{\bf Claim 9.} {\em Let $p>0$. For every $K\in\K^{n-1}$, and every $t_1,t_2\in\R$ with $t_1\le t_2$ we have:
\begin{eqnarray}\label{cylinders2}
\sigma_r(K\times[t_1,t_2])=\mu_r(w+I_{K\times [t_1,t_2]})=V_{n-1}(K)(\bar f(t_1,p)-\bar f(t_2,p)).
\end{eqnarray}}
\begin{proof}
We have, for every $t_1\in\R$,
$$
K\times [t_1,\infty)=\epi(u)\quad\mbox{where}\quad u=t_1+I_K.
$$
\begin{eqnarray}\label{cylinders1}
\mu_r(w+I_{K\times [t_1,\infty)})&=&\mu_r(w+I_{\epi(u)})\\
&=&\bar\mu(u)=\int_{\dom(u)}\bar f(u(x'),p)dx'\nonumber\\
&=&\int_K \bar f(t_1,p)dx'=V_{n-1}(K) f(t_1,p).\nonumber
\end{eqnarray}
On the other hand, for $t_1,t_2\in\R$ with $t_1\le t_2$ we have
$$
K\times [t_1,t_2]\cup K\times [t_2,\infty)=K \times[t_1,\infty),\quad
K\times [t_1,t_2]\cap K\times [t_2,\infty)=K\times \{t_2\};
$$
so that 
\begin{eqnarray*}
(w+I_{K\times [t_1,t_2]})\wedge(w+I_{K\times [t_2,\infty)})=w+I_{K\times [t_1,\infty)},\\
(w+I_{K\times [t_1,t_2]})\vee(w+I_{K\times [t_2,\infty)})=w+I_{K\times \{t_2\}}.
\end{eqnarray*}
Hence, as $\mu_r$ is a valuation and it is simple, and as $\dim(K\times \{t_2\})\le n-1$, we obtain
\begin{eqnarray*}
\mu_r(w+I_{K\times [t_1,t_2]})&=&\mu_r(w+I_{K\times [t_1,\infty)})-\mu_r(w+I_{K\times [t_2,\infty)})\\
&=&V_{n-1}(K)(\bar f(t_1,p)-\bar f(t_2,p)).\nonumber
\end{eqnarray*}
\end{proof}

The next step is to deduce further information about $\bar f$ exploiting the homogeneity of $\mu_r$ (recall
that $\mu_r$ is homogeneous of order $k$).

\bigskip

\noindent{\bf Claim 10.} {\em There exists a non-negative decreasing function $\phi\,:\,\R\to\R$ such that
$$
\bar f(t,p)=p^{n-1-k}\phi(tp)\quad\forall\,(t,p)\in\R\times(0,\infty).
$$}
\begin{proof}
We recall that $w_p(x)=p(x,e_n)$ for every choice of $p\ge0$ and $x\in\R^n$.
As before, let $K\in\K^{n-1}$ and let $\lambda>0$; we have, for $x\in\R^n$ and $t\in\R$,
\begin{eqnarray*}
w_p\left(\frac x\lambda\right)+I_{K\times [t,\infty)}\left(\frac x\lambda\right)&=&
w_{p/ \lambda}(x)+I_{\lambda (K\times [t,\infty))}(x)\\
&=&w_{p/ \lambda}(x)+I_{\lambda K\times [\lambda t,\infty)}(x).
\end{eqnarray*}

By the homogeneity of $\mu_r$
\begin{eqnarray*}
\mu_r\left((w_p+I_{K\times [t,\infty)})\left(\frac\cdot\lambda\right)\right)&=&
\lambda^{k}\mu_r(w_p+I_{K\times [t,\infty)})\\
&=&\lambda^{k} V_{n-1}(K) \bar f(t,p),
\end{eqnarray*}
and
\begin{eqnarray*}
\mu_r(w_{p/\lambda}+I_{\lambda K\times [\lambda t,\infty)})&=&
V_{n-1}(\lambda K)\bar f\left(\lambda t,\frac p\lambda\right)\\
&=&\lambda^{n-1}V_{n-1}(K)\bar f\left(\lambda t,\frac p\lambda\right)
\end{eqnarray*}
by the homogeneity of intrinsic volumes. Hence, as we may chose $K$ so that $V_{n-1}(K)>0$,
we obtain that for every $t\in\R$, $p>0$ and $\lambda>0$ we have
$$
\bar f(t,p)=\lambda^j \bar f\left(\lambda t,\frac p\lambda\right).
$$
with
$$
j=n-1-k.
$$
Taking $\lambda=p$ yields 
$$
\bar f(t,p)=p^j \bar f (tp,1)=p^j \phi(tp),
$$
where we have set 
$$
\phi(s)=\bar f(s,1)
$$
for all real $s$. As $\bar f$ is non-negative and decreasing with respect to $t$ for every $p>0$ the claim follows.
\end{proof}

In the next step we prove that the m-continuity of $\mu_r$ implies that the function $\phi$ is constant.

\bigskip

\noindent{\bf Claim 11.} {\em The function $\phi$ introduced in the previous step is constant in $\R$ (in particular $\sigma_r$ vanishes on cylinders).}

\begin{proof} By the previous steps we have that for every $K\in\K^{n-1}$ and for every $t_1,t_2\in\R$ with $t_1\le t_2$,
\begin{equation}\label{Claim 10}
\mu_r(w+I_{K\times[t_1,t_2]})=V_{n-1}(K)(\phi(pt_1)-\phi(pt_2))\,p^j.
\end{equation}
Let $K$ be the $(n-1)$-dimensional unit cube with centre at the origin and let
$$
D=\{x=(x_1,\dots,x_n)\,:\,(x_1,\dots,x_{n-1})\in K,\, -1\le x_n\le0\}.
$$
We also set, for $i\in\N$,
$$
E_i=D\cap\left\{x\,:\,-1\le x_n\le-\frac1i\right\},\quad
F_i=D\cap\left\{x\,:\,-\frac1i\le x_n\le0\right\}.
$$


\begin{figure}[h]
    \centering
    \includegraphics[width=0.5\textwidth]{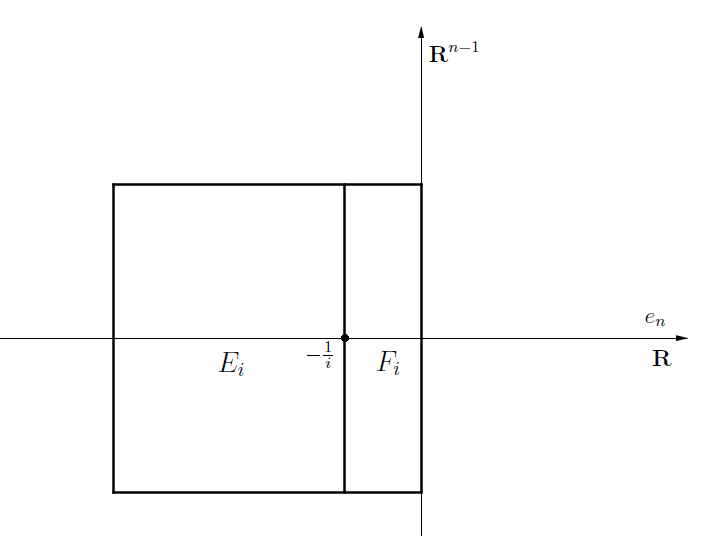}
    \caption{The sets $E_i$ and $F_i$}
    \label{}
\end{figure}

\bigskip

In particular, for every $i$,
$$
E_i,F_i\in\K^n,\quad E_i\cup F_i=D,\quad
\dim(E_i\cap F_i)=n-1.
$$
Let $s>0$. For $i\in\N$ define the function $v_i\,:\,\R^n\to\R$ as
$$
\bar v_i(x)=\bar v_i(x_1,\dots,x_n)=
s\cdot i\cdot\left(x_n+\frac1i\right)
$$  
and 
$$
v_i=\bar v_i\vee I_D.
$$
Note that 
$$
\mbox{$v_i=\infty$ in $\R^n\setminus D$,}\quad
\mbox{$v_i=0$ in $E_i$,}\quad
\mbox{$v_i=\bar v_i$ in $F_i$.}
$$
In particular $v_i$ is a decreasing sequence of functions in $\C^n$ converging to $I_D$ in the relative interior of $D$, so that by m-continuity we have
$$
\lim_{i\to\infty}\mu_r(v_i)=\mu_r(I_D)=0,
$$
where we have used the fact that $\mu_r$ vanishes horizontally (Claim 1). We may also write
$$
v_i=(\bar v_i+I_{F_i})\vee I_{E_i},
$$
and using the fact that $\mu_r$ is a simple valuation, and Claim 1 again, we get that $\mu_r(v_i)=
\mu_r(\bar v_i+I_{F_i})$ so that
$$
\lim_{i\to\infty}\mu_r(\bar v_i+I_{F_i})=0.
$$
On the other hand, by translation invariance, if we set 
$$
u_i(x)=u_i(x_1,\dots,x_n)=
\left(\bar v_i+I_{F_i}\right)\left(x_1,\dots,x_{n-1},x_n-\frac1i\right),
$$
we find that 
$$
u_i=w_i+I_{K\times[0,t_i]},
$$
where
$$
w_i(x)=si\,x_n\quad\mbox{and}\quad
t_i=\frac1i.
$$
Consequently, by \eqref{Claim 10}
$$
\mu_r(u_i)=(si)^j\, V_{n-1}(K)(\phi(0)-\phi(s))\quad\forall\, i\in\N.
$$


\begin{figure}[h]
    \centering
    \includegraphics[width=0.48\textwidth]{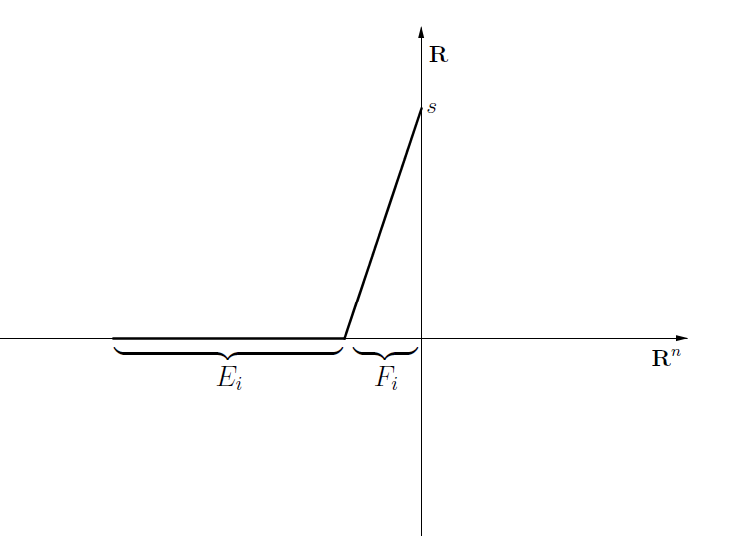}
    \includegraphics[width=0.48\textwidth]{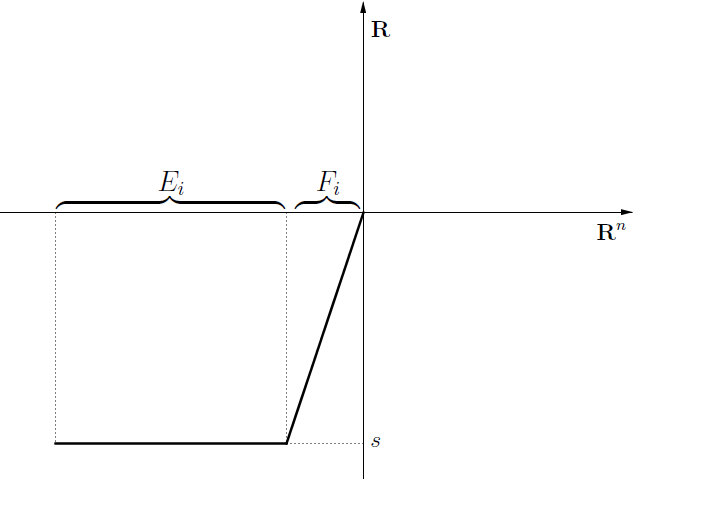}	
    \caption{The construction of $v_i$ for $s>0$ and $s<0$}
    \label{}
\end{figure}

Letting $i$ tend to infinity this quantity must tend to zero, by the previous part of the proof; as $j\ge0$ and $V(K)>0$, the only 
possibility is $\phi(s)=\phi(0)$. This proves that $\phi$ is constantly equal to $\phi(0)$ in $[0,\infty)$. 

To achieve the same result in $(-\infty,0]$ we may argue in a similar way. Let $s<0$ and $K$, $D$, $E_i$ as above. Set
$$
\bar v_i(x_1,\dots,x_n)=-si\,x_n
$$ 
and 
$$
v_i=(I_D+s)\vee \bar v_i.
$$
This is again a decreasing sequence in $\C^n$, converging to $s+I_D$ in the relative interior of $D$; by Claim 1:
$$
\lim_{i\to\infty}\mu_r(v_i)=0.
$$
On the other hand
$$
\mu_r(v_i)=\mu_r(\bar v_i+I_{K\times[-1/i,0]})
=(-si)^jV_{n-1}(K)(\phi(s)-\phi(0)).
$$
The conclusion $\phi(s)=\phi(0)$ follows as above.

\end{proof}

\bigskip

\noindent{\bf Claim 12.} {\em Let $V\in\R^n$, $c\in\R$ and $K\in\K^n$; define
$$
u\,:\,\R^n\to\R,\quad
u(x)=(x,V)+c+I_K(x).
$$
Then
$$
\mu_r(u)=0.
$$}

\begin{proof}
Assume first that $c=0$.
If $V=0$ the assert follows from Claim 1. Assume $V\ne0$ and let $V=pe$, with $p>0$ and $e$ a unit vector. Recalling the definition of 
$\sigma_{\mu_r,V}$ we have:
$$
\mu_r(u)=\sigma_{\mu_r,V}(K).
$$
On the other hand, since $\mu_r$ is rigid motion invariant, we can assume, without loss of generality, that $V=e_n$ and, as remarked in Claim 11, setting as before $w\,:\,\R^n\to\R$ defined by $w(x)=(x,pe_n)$, we get
$$
\sigma_{\mu_r,V}(H\times [t_1,t_2])=
\mu_r(w+I_{H\times [t_1,t_2]})=0.
$$
for every $H\in\K^{n-1}$, $t_1,t_2\in\R$ such that $t_1\le t_2$. Let us choose $H$, $t_1$ and $t_2$ such that
$$
K\subseteq H\times[t_1,t_2].
$$
Then, as $\sigma_{\mu_r,V}$ is non-negative and monotone increasing (Claims 4 and 5),
$$
0\le\sigma_{\mu_r,V}(K)\le\sigma_{\mu_r,V}(H\times [t_1,t_2])=0.
$$
The case $c\ne0$ is readily recovered by the previous one using the translation invariance of $\mu_r$.
\end{proof}

The last result will open the way to prove that $\mu_r$ vanishes on piece-wise linear functions and, eventually, it vanishes 
identically on $\C^n$. 

\begin{definition} A function $u\in\C^n$ is said to be piece-wise linear if:
\begin{itemize}
\item $\dom(u)=P$ is a polytope;
\item there exists a polytopal partition
${\mathcal P}=\{P_1,\dots,P_N\}$ of $P$ such that for every $i\in\{1,\dots,N\}$ there exists $V_i\in\R^n$ and $c_i\in\R$ such that
$$
u(x)=c_i+(x,V_i)\quad\forall\,x\in P_i.
$$
\end{itemize}

\end{definition}

\bigskip

\noindent{\bf Claim 13.} {\em The valuation $\mu_r$ vanishes on piece-wise linear functions.} 

\begin{proof} As any polytopal partition admits a refinement which is a complete partition (see Remark \ref{complete refinement} in section \ref{subsection Complete partitions}), 
without loss of generality we may assume that $\mathcal P$ is complete, so that in particular it is an inductive partition (see Proposition
\ref{complete partitions are D-partitions}). The claim follows immediately from Claim 12, the fact that $\mu_r$ is simple, and 
Lemma \ref{decomposition of simple valuations}.
\end{proof}

\bigskip

\noindent{\bf Claim 14.} {\em The valuation $\mu_r$ vanishes on $\C^n$.} 

\begin{proof} Let $u\in\C^n$; if the dimension of $\dom(u)$ is strictly less than $n$, $\mu_r(u)=0$ as $\mu_r$ is simple. So, assume that
$\Omega=\interno(\dom(u))\ne\emptyset$. Let $P$ be a polytope contained in $\Omega$, and let $u_i$, $i\in\N$, be a sequence 
of piece-wise linear functions of $\C^n$, such that for every $i$: $\dom(u_i)=P$, $u_i\ge u_{i+1}$ in $P$, and the sequence $u_i$ 
converges uniformly to $u$ in $P$; such a sequence exists by standard approximation results of convex functions by piece-wise linear functions.
Using the m-continuity of $\mu_r$ and the previous Claim 13, we obtain
$$
\mu_r(u+I_P)=\lim_{i\to\infty} \mu_r(u_i)=0.
$$
Now take a sequence of polytopes $P_i$, $i\in\N$, such that: $P_i\subseteq P_{i+1}\subseteq\Omega$ for every $i$ and
$$
\Omega=\bigcup_{i\in\N} P_i.
$$
Then the sequence
$$
u+I_{P_i},\quad i\in\N,
$$
is formed by elements of $\C^n$, is decreasing, and converges point-wise to $u$ in $\Omega$; by m-continuity and the previous part of this proof
$$
\mu_r(u)=\lim_{i\to\infty}\mu_r(u+I_{P_i})=0.
$$
\end{proof}

The proof of Theorem \ref{characterization of k-homogeneous valuations} is complete, under the additional assumption that the density $f$ of $\mu$ is smooth. 
The next and final step explains how to deduce the theorem in the general case.

\bigskip

\noindent{\bf Claim 15.} {\em The assumption that $f$ is smooth can be removed.} 

\begin{proof} Let $\mu$ be as in the statement of Theorem \ref{characterization of k-homogeneous valuations}, and let $\mu_i$, $i\in\N$, be the sequence of valuations
 determined by Proposition \ref{mollification} (taking for example $\epsilon=1/i$, $i\in\N$). It follows from the definition of $\mu_i$ given in section \ref{regularization of geometric densities} that, as $\mu$ is 
 $k$-homogeneous, $\mu_i$ is $k$-homogeneous as well. Moreover, the only  non-vanishing geometric density of $\mu_i$, that we will denote by $f_i$, is smooth.
 Hence, for every $i$ we may apply the previous part of the proof to $\mu_i$ and deduce that
$$
\mu_i(u)=\int_{\R}V_k(\cl(\{u< t\}))d\nu_i(t),
$$
where $\nu_i$ is a Radon measure on $\R$ and it is related to $f_i$ by the equality
$$
f_i(t)=\int_{(t,\infty)} d\nu_i(s),\quad\forall\, t\in\R.
$$
We apply Proposition \ref{Wright representation} to get
$$
\mu_i(u)=\int_{\R}f_i(t)d\beta_k(u;t)\quad\forall\,i\in{\N},\, u\in\C^n,
$$
we recall that $\beta_k(u;\cdot)$ is the distributional derivative of the increasing function
$$
\R\ni t\,\mapsto\, V_k(\cl(\{u< t\})).
$$
From Proposition \ref{structure beta_i} we know that $\beta_k(u;\cdot)$
can be decomposed as the sum of a part which is absolutely continuous with respect to the one-dimensional Lebesgue 
measure and a Dirac point-mass measure having support at $m(u)$ and weight $V_k(\{x\,:\,u(x)=m(u)\})$. If in particular we assume that $u\in\C^n$ is such that
\begin{equation}\label{claim 14}
\{x\,:\,u(x)=m(u)\}\;\mbox{consists of a single point,}
\end{equation}
we have that (as $k\ge1$) 
$$
V_k(\{u=m(u)\})=0,
$$
so that $\beta_k(u;\cdot)$ is absolutely continuous with respect to the Lebesgue measure on the real line.

Our next move is to prove that, under the assumption \eqref{claim 14}
\begin{equation}\label{claim 14.2}
\lim_{i\to\infty}\int_{\R}f_i(t)d\beta_k(u;t)=
\int_{\R}f(t)d\beta_k(u;t).
\end{equation}
We know that the sequence $f_i$ converges to $f$ almost everywhere on $\R$ with respect to the Lebesgue measure, and hence with respect to
$\beta_k(u;\cdot)$. Note also that
$$
f_i(t)=\int_\R f(t-s)g_{1/i}(s)ds=
\int_{-1}^1 f(t-s)g_{1/i}(s)ds
$$
where $g$ is the mollifying kernel introduced in section \ref{regularization of geometric densities} (which in particular is supported in $[-1,1]$)
and
$$
g_\epsilon(s)=\frac1\epsilon g\left(\frac s\epsilon\right),\quad\forall\,\epsilon>0.
$$ 
As $f$ is decreasing (and non-negative)
$$
0\le f_i(t)\le
\int_\R f(t-1)g_{1/i}(s)ds=f(t-1)\quad\forall
t\in\R,\,i\in\N.
$$
On the other hand
\begin{eqnarray*}
\int_{\R}f(t-1)d\beta_k(u;t)&=&
\int_\R f(t)d\beta_k(u;t+1)\\
&=&\int_\R f(t)d\beta_k(\bar u;t)=
\int_\R V_k(\cl(\{\bar u< t\}))d\nu(t)<\infty,
\end{eqnarray*}
where $\bar u=u-1$ and the last inequality is due to the integrability condition on $f$ (Proposition \ref{integrability condition} and Corollary \ref{corollary k-hom}). 
Hence we may apply the dominated convergence theorem and obtain \eqref{claim 14.2}. Note that if $u$ verifies condition \eqref{claim 14},
then so does the function $u+s$, for every $s\in\R$. By Proposition \ref{mollification} we conclude that 
\begin{equation}\label{claim 14.3}
\mu(u+s)=\int_\R f(t)d\beta_k(u+s;t)=
\int_\R f(t)d\beta_k(u;t-s)=
\int_\R f(t+s)d\beta_k(u;t),\quad\mbox{for a.e. $s\in\R$.} 
\end{equation}
Let $s_i$, $i\in\N$, be a decreasing sequence of real numbers converging to zero such that \eqref{claim 14.3} holds true; then by 
$m$ continuity of $\mu$
$$
\lim_{i\to\infty}\mu(u+s_i)=\mu(u).
$$
The m-continuity implies also that $f$ is right-continuous (see right after Corollary \ref{m-continuita' e continuita' destra}), hence
$$
\lim_{i\to\infty} f(t+s_i)=f(t)\quad\forall\, t\in\R.
$$
Using again the monotonicity of $f$ and the monotone convergence theorem we obtain
$$
\lim_{i\to\infty}\int_\R f(t+s_i)d\beta_k(u;t)
=\int_\R f(t)d\beta_k(u;t).
$$
Putting the last equalities together we arrive to
\begin{equation}\label{claim 14.4}
\mu(u)=\int_\R f(t)d\beta_k(u;t)
=\int_\R V_k(\cl(\{u< t\}))d\nu(t),
\end{equation}
for every $u\in\C^n$ verifying \eqref{claim 14}. The last step will be to prove that this equality is true for every $u\in\C^n$. For 
$i\in\N$ set
$$
u_i\,:\,\R^n\to\,\R\quad
u_i(x)=u(x)+\frac{|x|^2}{i}.
$$ 
Clearly $u_i\in\C^n$ and, as $u_i$ is strictly convex it verifies condition \eqref{claim 14} and, consequently, \eqref{claim 14.4}. By m-continuity
$$
\lim_{i\to\infty}\mu(u_i)=\mu(u).
$$
We need to prove that
\begin{equation}\label{Claim 15 100}
\lim_{i\to\infty}\int_\R V_k(\cl(\{u_i< t\}))d\nu(t)=
\int_\R V_k(\cl(\{u<t\}))d\nu(t).
\end{equation}
As $u_i\ge u$ in $\R^n$ fo every $i$ we have that $\{u_i< t\}\sub\{u< t\}$ for every $t$.
We have already proven that
$$
\lim_{i\to\infty}\cl(\{u_i<t\})=\cl(\{u<t\})\quad\forall\, t\in\R,
$$
where the limit is intended in the Hausdorff metric on $\K^n$. 
Then \eqref{Claim 15 100} follows by the monotone the monotone convergence theorem,
and Theorem \ref{characterization of k-homogeneous valuations} is finally proven in the general case as well.
\end{proof}

\section{A non level-based valuation}
In this section we will present a way to construct monotone valuations on $\C^n$ which are moreover rigid motion invariant and m-continuous and, despite verifying all these desirable properties, cannot be expressed as a linear combination of homogeneous valuations on $\C^n$. 

Fix $n,m\in\N$, for all $u\in\C^n$ we set $\hat{u}(x,y)=u(x)+|y|$ for all $(x,y)\in\R^n\times \R^m$, where $|\cdot|$  is to be interpreted as the Euclidean norm in $\R^m$. Note that if $u\in\C^n$, then $\hat{u}\in\C^{n+m}$.

We are now ready to define the prototype of the valuations described at the beginning of this section.

\begin{proposition} \label{valbrutt} 
Let $n,m$ be fixed natural numbers and let $k\in\{0,\dots,n+m\}$. Let $t\in\R$. 
Then the map $\C^n\to \R$, defined as $u\mapsto V_k(\cl(\{\hat{u}<t\}))$, 
\begin{itemize}
\item[{\em i)}] is a valuation,
\item[{\em ii)}] is monotone decreasing,
\item[{\em iii)}] is rigid motion invariant,
\item[{\em iv)}] is m-continuous.
\end{itemize}
\end{proposition}
\begin{proof}
First of all, for ease of notation, set $\mu(\cdot)=V_k(\cl(\{\cdot <t\}))$.
By Proposition \ref{atomo}, $\mu$ verifies all the properties {\em i)} - {\em iv)}.

{\em i)} As a preliminary step to prove the condition on $\boldsymbol{\infty}$, notice that $\widehat{\boldsymbol{\infty}}=\boldsymbol{\infty}\in\C^{n+m}$. As a matter of fact, for all $(x,y)\in\R^n\times \R^m$ we have
$$\widehat{\boldsymbol{\infty}}(x,y)=\boldsymbol{\infty}(x)+|y|=\infty+|y|=\infty. $$
As a consequence 
$$ \mu(\widehat{\boldsymbol{\infty}})=\mu(\boldsymbol{\infty})=0. $$
Let now $u,v\in\C^n$, we have
\begin{subequations}
\begin{align}
\widehat{u\land v}=\hat{u}\land \hat{v}, \label{prima} \\
\widehat{u\lor v}=\hat{u}\lor \hat{v}. \label{seconda} 
\end{align}
\end{subequations}
We are going to prove \eqref{prima} only, as \eqref{seconda} is completely analogous.
Let $u,v\in\C^n$, then, for all $(x,y)\in\R^n\times \R^m$ we get
$$
\hat{u}(x,y)\land \hat{v}(x,y)= \left(u(x)+|y|\right)\land \left(v(x)+|y|\right)= 
u(x)\land v(x) + |y|= \widehat{u\land v}(x,y),
$$
and \eqref{prima} is proven.
Let now $u,v\in\C^n$ with $u\land v \in \C^n$:
$$
\mu(\widehat{u\lor v})+\mu(\widehat{u\land v})=\mu(\hat{u}\lor\hat{v})+\mu(\hat{u}\land \hat{v})=\mu(\hat{u})+\mu(\hat{v}),
$$
where in the last equality we have employed the valuation property of $\mu$.

{\em ii)} Let $u,v\in\C^n$ with $u\le v$.
It is immediate to verify that $\hat{u}\le \hat{v}$. Indeed, for all $(x,y)\in\R^n\times \R^m$,
$$ \hat{u}(x,y)=u(x)+|y|\le v(x)+|y|=\hat{v}(x,y).$$ 
As $\mu$ is monotone we have
$$u\le v \implies \hat{u}\le \hat{v}\implies \mu(\hat{u})\ge \mu(\hat{v}).$$

{\em iii)} Let $T$ be a rigid motion of $\R^n$ and let $u$ be a convex function in $\C^n$.
We define $u_T(x)=u(T(x))$ for all $x\in\R^n$.
We have 
$$\widehat{u_T}(x,y)=u_T(x)+|y|=u(T(x))+|y|=\hat{u}_{\hat{T}}(x,y)$$ 
for all $(x,y)\in\R^n\times \R^m$ where $\hat{T}$ is the rigid motion of $\R^n\times \R^m$ defined by $(x,y)\mapsto (T(x),y)$.
Therefore, as $\mu$ is rigid motion invariant,
$$\mu(\widehat{u_T})=\mu(\hat{u}_{\hat{T}})=\mu(\hat{u}).$$
In other words, the map $u\mapsto \mu(\hat{u})$ is rigid motion invariant as claimed.

{\em iv)} Let $u\in\C^n$ and let $u_i$, $i\in\N$, be a point-wise decreasing sequence of convex functions in $\C^n$ converging to $u$ point-wise in $\relint(\dom (u))$. We 
want to show that $\lim_{i\to\infty}\mu(\widehat{u}_i)=\mu(\hat{u})$.
In order to prove it, we will use the m-continuity of $\mu$ and show that the sequence $\widehat{u}_i$, $i\in\N$, is also a point-wise decreasing sequence of convex 
functions (this time in $\C^{n+m}$) that converges to $\hat{u}$ in the relative interior of its domain. By the reasoning used to prove {\em ii)} we deduce that $\widehat{u}_i$, 
$i\in\N$, is point-wise decreasing as well.
Notice also that 
\begin{equation}
\begin{split}
\dom(\hat{u})=\{(x,y)\in\R^n\times \R^m \,:\, u(x)+|y|<\infty \} \\
= \{(x,y)\in\R^n\times \R^m \,:\, u(x)<\infty\}=\dom(u)\times \R^m;
\end{split}\nonumber
\end{equation}
so that $$\relint(\dom(\hat{u}))=\relint(\dom(u))\times \R^m.$$
Let now $(x,y)\in\relint(\dom(u))$, then  
$$ \widehat{u_i}(x,y)=u_i(x)+|y|\to u(x)+|y|=\hat{u}(x,y).$$
We conclude by the m-continuity of $\mu$.
\end{proof}

Let us specialize the valuation of Proposition \ref{valbrutt} to the case $n=m=k=1$.
In this case we can provide a simple geometric explanation: $V_1(\cl(\{\hat{u}<t\}))$ is equal to the length of that portion of the graph of $u$ that lies strictly under the level $t$ and therefore we will refer to it as \emph{undergraph-length}.

To see that, first consider $t\le m(u)=m(\hat{u})$: in this case the set $\{\hat{u}<t\}$ is empty and so $V_1(\cl(\{\hat{u}<t\}))$ trivially equals the length of the graph lying strictly under $t$, the latter being $0$ as well.
On the other hand, let $t>m(u)=m(\hat{u})$; then, by Corollary \ref{corollario sottolivelli} we have that 
$\cl(\{\hat{u}<t\})=\{\hat{u}\le t\}$. Note that this set can be rewritten as 
$$\{ (x,y)\in\R^2 \,:\, u(x)+|y|\le t\}=\{(x,y)\in\R^2 \,:\, |y|\le t-u(x)\}.$$
In other words, $\cl(\{\hat{u}<t\})$ be obtained as a result of the following process: take the part of $\epi(u)$ that lies below the line $\{(x,y)\in\R^2 \,:\, y= t\}$, translate it 
``vertically'' so that the flat top is now lying on the $x$-axis $H:=\{(x,0)\in\R^2\}$, finally symmetrize it with respect to $H$. We recall that $V_1(K)$ coincides with the length 
(1-dimensional Lebesgue-measure) in case $\dim(K)=1$ and with $\frac{1}{2} \mathcal{H}^1(\partial K)$ when $\dim(K)=2$ (see \cite{Schneider}).
If $\dom(u)$ has dimension 1, as $t>m(u)$, the $\epi(u)\cap \{(x,y)\in\R^2 \,:\, y\le t\}$ is 2-dimensional and therefore $V_1(\cl(\{\hat{u}<t\}))$ is equal to 
the length of the graph of $u$ that lies strictly under the level $t$ for every choice of $t\in\R$.

The undergraph-length is not a level based valuation. By these words we mean that we could actually take a convex function $u\in\C^1$, rearrange its levels using translations 
and obtain another convex function $v$ such that $V_1(\cl(\{\hat{u}<t\}))\ne V_1(\cl(\{\hat{v}<t\}))  $ for all $t>m(u)$.   
Take for instance $u(x)=|x|$ and $v(x)= x/2+I_{[0,\infty)}$ for all $x\in\R$; we have $\{u<t\}=(-t,t)$ and $\{v<t\}=(0,2t)=t+(-t,t)$ for all positive real $t$. On the other hand, their 
undergraph-lengths differ: a quick use of the Pythagorean theorem reveals that $V_1(\cl(\{\hat{u}<t\}))=2\sqrt{2}t$ while  $V_1(\cl(\{\hat{v}<t\}))=\sqrt{5}t $.  

\par The length of the undergraph is a valuation which is completely different from the ones we have studied so far: not only it is not $\alpha$-homogeneous for any real $
\alpha$, it turns out that $V_1(\cl(\{\,\widehat{\cdot}<t\})) $ cannot even be written as a finite sum of homogeneous functions.
To prove this consider the following $u\in\C^1$, defined as $u(x)=|x|$ for all $x\in\R$. 
For all $\lambda >0 $ we get $V_1(\cl(\{\widehat{u_\lambda}<1\}))=2\sqrt{1+\lambda^2} $.
Since $V_1(\cl(\{\widehat{u_\lambda}<1\}))$ is not a polynomial in $\lambda$,  $V_1(\cl(\{\widehat{u_\lambda}<1\}))$ cannot be decomposed into the
(finite) sum of homogeneous 
functions. This implicitly tells us that under these assumptions (monotonicity, rigid motion invariance and m-continuity), homogeneous valuations do not form a basis for the 
vector space of valuations on $\cn$.

\bigskip

\noindent
L. Cavallina and A. Colesanti:\ Dipartimento di Matematica e Informatica ``U.Dini", Viale Morgagni 67/A, 50134, Firenze, Italy

\noindent
Electronic mail addresses:
cavathebest@hotmail.it, colesant@math.unifi.it

\end{document}